\documentclass[11pt]{article}
\usepackage{graphicx}
\usepackage{amssymb}
\usepackage{epstopdf}
\input xy
\xyoption{all}

\newtheorem{lemma}{Lemma}[section] 
\newtheorem{propos}[lemma]{Proposition}
\newtheorem{example}[lemma]{Example}
\newtheorem{theorem}[lemma]{Theorem}
\newtheorem{cor}[lemma]{Corollary}

\newtheorem{defin}[lemma]{Definition}

\newcommand{\C}{\mathbb{C}}

\newcommand{\F}{\mathbb{F}}

\newcommand{\proof}{{\noindent {\bfseries  Proof:}\quad }}

\newcommand{\ev}{\mathrm{ev}}
\newcommand{\coev}{\mathrm{coev}}

\newcommand{\extd}{\mathrm{d}}

\newcommand{\eps}{{\epsilon}}
\newcommand{\tens}{\mathop{\otimes}}
\newcommand{\la}{{\triangleright}}
\newcommand{\ra}{{\triangleleft}}

\newcommand{\id}{\mathrm{id}}

\newcommand{\<}{\langle}
\renewcommand{\>}{\rangle}
\newcommand{\End}{\mathrm{ End}}

\begin{document}
\title{Line bundles and the Thom construction in noncommutative geometry}
\author{E.J.\ Beggs\footnote{E-mail: E.J.Beggs@swansea.ac.uk}~ \& Tomasz Brzezi\'nski\footnote{E-mail: T.Brzezinski@swansea.ac.uk} \\
Department of Mathematics, Swansea University, Swansea SA2 8PP, U.K.}

\maketitle

\begin{abstract}  The idea of a line bundle in classical geometry is transferred to noncommutative geometry by the idea of a Morita context. From this we can construct  $\mathbb{Z}$- and $\mathbb{N}$-graded algebras, the $\mathbb{Z}$-graded algebra being a Hopf-Galois extension. 
A non-degenerate Hermitian metric gives a star structure on this algebra, and an additional star operation on the line bundle gives a star operation on the $\mathbb{N}$-graded algebra. In this case, we can carry out the associated circle bundle and Thom constructions. Starting with a $C^*$-algebra as base, and with some positivity assumptions, the associated circle and Thom algebras are also $C^*$-algebras. We conclude by examining covariant derivatives and Chern classes on line bundles after the method of Kobayashi \& Nomizu.
\end{abstract}

\section{Introduction} In this paper we consider the generalisation of line bundles to noncommutative geometry. Of course, it has been successfully assumed for many years that finitely generated projective modules over an algebra $A$  should be understood as (sections of) vector bundles in noncommutative geometry. In Section \ref{bhakvck} we will characterise the specialisation of this idea which generalises a line, rather than a higher dimensional, bundle. To many readers much of the material in this section will seem familiar, and this will be explained in Section \ref{vcuadrts}, where we note that the definition of a line module is just that of a Morita context, and that the results of Section \ref{bhakvck} are essentially taken out of the literature on algebraic $K$-theory. Our order of exposition has been chosen to make things slightly clearer to readers who are not familiar with Morita theory. 

At this point, the reader may ask why we restrict a perfectly satisfactory theory to a special case, especially if that special case has already been studied (though in a rather different setting). The answer is that we can construct two algebras of functions from the line module, an $\mathbb{N}$-graded and a $\mathbb{Z}$-graded
algebra. The $\mathbb{Z}$-graded algebra can be characterised using the theory of Hopf-Galois extensions, culminating in Theorem \ref{vchuaytdrtys}, which gives
a 1-1 correspondence between automorphisms of the category ${}_A\mathbb{M}$ of left $A$-modules, left line modules over $A$, and Hopf-Galois $\mathbb{CZ}$ extensions of $A$. 

To further study these algebras, we introduce a Hermitian metric, or inner product on the module, which is essentially given by a Hilbert $C^*$-module.
(See \cite{BegMa4} for more of this approach, and \cite{Lance} for Hilbert $C^*$-modules.) The Hermitian metric is described in Section~\ref{gahcvjdvyu}, and is used to
give a star operation on the $\mathbb{Z}$-graded algebra in Section~\ref{vuwicvhcv}. To give a star operation on the $\mathbb{N}$-graded algebra requires another structure, an involution on the module itself, and this is described in Section~\ref{xctkvdft}. The $\mathbb{N}$-graded algebra 
does not have a direct interpretation as a Hopf-Galois extension, but with the additional structure of a Hermitian metric and a star operation, we can form a 
 $\mathbb{Z}/2$ Hopf-Galois extension from it, see Section~\ref{cvaufv2}.

These two algebras (and their different star operations) correspond to two different cases in the classical theory. The complex canonical line bundle in algebraic geometry has a star operation, given by conjugation of the complex coordinates, but that operation takes values in the dual of the canonical bundle rather than the original bundle. 
This corresponds to the $\mathbb{Z}$-graded algebra and involution. However the case of a real line bundle on a topological space is rather different. The functions
on the total space of the bundle which are polynomial in the fibre $\mathbb{R}$ form an $\mathbb{N}$-graded algebra. The involution is just complex conjugation of functions on the total space of the bundle. This corresponds to the $\mathbb{N}$-graded algebra and involution.

Now we come to a case of the Thom construction \cite{Thom}. Given an $\mathbb{R}^n$ bundle on a compact topological space, there are two obvious compactifications to consider. One (the Thom space) is the one point compactification of the total space, and the other adds a point to each fibre to make an associated sphere $S^n$ bundle. We shall consider the case of a line bundle, $n=1$. The $\mathbb{N}$-graded algebra of functions polynomial in the fibre direction can be modified to 
give an algebra of functions vanishing at infinity on the fibres. Making this into a unital algebra by adding $\mathbb{C}$ corresponds to the Thom construction, whereas adding a copy of $A$ corresponds to the associated sphere (or circle in this case) bundle. 

Another successful assumption about noncommutative geometry is that compact Hausdorff topological (noncommutative) spaces correspond to unital $C^*$-algebras. If the original algebra $A$ is a unital $C^*$-algebra, we might hope that the algebras corresponding to the Thom construction and associated circle bundle might also be unital $C^*$-algebras. In Sections~\ref{cvaufv1} and \ref{cvaufv3} we show that this is indeed the case, by explicit construction of a star representation of the $\mathbb{N}$-graded algebra (or rather, its vanishing at infinity variant) on a Hilbert space. 

 In Section~\ref{cgfhjhfxsz} we discuss automorphisms of the centre of the algebra, and how they arise from line modules. This was extensively studied in  \cite{FroPicard}. 
 This is used in various other parts of the paper, as central elements give the freedom of choice in constructing bimodule maps between line modules. Section~\ref{pichgcjf} on the Picard group of an algebra (see \cite{BassK}) contains more questions than answers, but hopefully they may be useful questions! The Picard group consists of isomorphism classes of line modules. The basic problem is whether there exists a `universal' Hopf-Galois extension of the algebra $A$ which gives all line modules over $A$. In particular the problem of constructing a Hopf-Galois extension with group the Picard group (or some related extension) may be bound up in choosing cocycles with values in the centre. 

In Section~\ref{constcxagjy} we look at a method of constructing line modules which has arisen in the literature (e.g.\ \cite{LandiDecon,Landi4DHopf}), and give a simple example, the noncommutative Hopf fibration. Finally in Section~\ref{acjfjdssd} we will look at the differential geometry of covariant derivatives on line modules. It turns out that we can give very explicit descriptions of these covariant derivatives. One reason to be interested in such things is that 
covariant derivatives on line bundles give $U(1)$ gauge theory, otherwise known as electromagnetism, in theoretical physics. We apologise to the reader in not being experts on the theoretical physics literature here, but \cite{MaScScWe} seems to be a good reference for the physics, and \cite{BrzMaj:gau} for the mathematics. The idea of bimodule covariant derivative mentioned here has a long history, see for example \cite{DVMic,DVMass,MouLC,MadIntro}. However in this article we concentrate on another application,
 following the example of classical differential geometry in \cite{KobNom}, we give a definition of the differential Chern class in terms of the trace of the curvature. 

For related material in the literature, in \cite{ConnesThom} there is a version of the Thom map for $K$-theory associated to a $C^*$ dynamical $\mathbb{R}$-system. This refers to the topological space results for the Thom isomorphism in $K$-theory given in \cite{AtiK}. For our current paper, it is of interest to note that the discussion in \cite{AtiK} makes use of {\em decomposable vector bundles}, i.e.\ vector bundles which are a sum of line bundles.

So what is the chance of extending these results to the equivalent of $\mathbb{R}^n$ bundles rather than just line bundles? The basic difficulty is what reduces to the relations between the `fibre coordinates'. In other words, the tensor algebras are far too big, and we have to perform what would be classically a symmetrisation operation. This may involve a map $\sigma:E\tens_A E\to E\tens_A E$, as may be indicated from the idea of bimodule covariant derivatives, or braidings associated to Yetter-Drinfeld modules (see \cite{worondiff}), but the idea of what would work in a large number of cases is not really understood.

One idea possibly worthy of future investigation is the relation to cohomology of algebras. For topological spaces, there is a 1-1 correspondence between principal abelian group bundles and 
$H^1$ cohomology with coefficients in the abelian group. It is extremely tempting to transfer this idea to algebras by taking the Hopf-Galois extensions in place of the principal abelian group bundles. The big question, of course, is whether this is worthwhile in terms of results, extension to $H^n$, and examples. The evidence presented here on $\mathbb{Z}$ and $\mathbb{Z}/2$ Hopf-Galois extensions being related to line modules and line modules with real structure, respectively, might be taken as supporting this idea.

The authors would like to thank David Evans, Jeffrey Giansiracusa, Ulrich Kr\"ahmer, Giovanni Landi  and Claudia Pinzari for their help during the preparation of this paper.

\section{Prerequisites}
We work over the field $\mathbb{C}$. 
Begin with a unital algebra $A$, and suppose that $E$ is a left $A$-module. 

\begin{defin} \label{vcgaftydgs}
 The \emph{(left) dual} $E^\circ$ of a left
$A$-module $E$ is defined to be ${}_AHom(E,A)$, the left module maps from $E$ to $A$.
Then $E^\circ$ has a right module structure given by $(\alpha.a)(e)=\alpha(e).a$ for all
 $\alpha\in E^\circ$, $a\in A$ and $e\in E$. 
\end{defin}

\begin{defin}\label{canonyy} A left $A$-module $E$ is said to be \emph{finitely 
    generated projective} if
there are $e^i\in E$ and $e_i\in E^\circ$ (for integer $1\le i\le n$)
(the `dual basis') so that for all $f\in E$, $f=\sum e_i(f).e^i$. 
From this it follows directly that $\alpha=\sum e_i.\alpha(e^i)$
for all $\alpha\in E^\circ$. The $A$ valued matrix $P_{qj}=\ev(e^q\tens e_j)$
is an idempotent associated to the module.
\end{defin}

If $E$ is a bimodule, then $E^\circ$ is also a bimodule, with right module structure as in 
Definition \ref{vcgaftydgs} and left module structure $(a.\alpha)(e)=\alpha(e.a)$. If in addition $A$ is
finitely generated projective as a left $A$-module, we have bimodule maps
\begin{eqnarray*}
\ev:E\tens_A E^\circ\to A\ ,\quad \coev:A\to E^\circ\tens_A E\ ,
\end{eqnarray*}
given by $\ev(e\tens\alpha)=\alpha(e)$ and $\coev(1_A)=\sum_i e_i\tens e^i$. These obey the identities
\begin{eqnarray} \label{xerakjcgh}
(\ev\tens\id)(\id\tens\coev(1_A))\ =\ \id:E &\to & E\ ,\cr
(\id\tens\ev)(\coev(1_A)\tens\id)\ =\ \id:E^\circ &\to & E^\circ\ .
\end{eqnarray}

Now suppose that $A$ is a star algebra. 
 We shall find it convenient to use the notation of conjugate modules, as given in \cite{BegMa2}. If $E$ is an $A$-bimodule, then the conjugate module $\overline{E}$ is an $A$-bimodule, with, for $a\in A$ and $\overline{e}\in \overline{E}$,
\begin{eqnarray*}
a.\overline{e}\ =\ \overline{e.a^*}\ ,\quad \overline{e}.a\ =\ \overline{a^*.e}\ .
\end{eqnarray*}
This notation has the benefit that many operations which do not look initially like bimodule maps can be written as bimodule maps, for example the star operation
on a star algebra $A$ can be written as a bimodule map $\star:A\to \overline{A}$,
$a\mapsto \overline{a^*}$. We also use the bimodule map, for $A$-bimodules $E$ and $F$, $\Upsilon: \overline{F\tens_A E} \to \overline{E}\tens_A \overline{F}$, with the inverse
\begin{eqnarray*}
\Upsilon^{-1}:\overline{E}\tens_A \overline{F}\to \overline{F\tens_A E}\ ,\quad
\Upsilon^{-1}(\overline{e}\tens_A \overline{f})\ =\ \overline{f\tens e}\ ,
\end{eqnarray*}
and the isomorphism $\mathrm{bb}:E\to \overline{\overline{E}}$
given by $\mathrm{bb}(e)=\overline{\overline{e}}$.

Given a differential calculus $(\Omega^nA,\extd,\wedge)$ defined in terms of a differential graded algebra with $\Omega^0 A=A$, we can define an $A$-covariant derivative on a left $A$-module $E$
as follows:

\begin{defin} Given a left $A$-module $E$, a {\em left $A$-covariant derivative}
is a map $\nabla:E\to\Omega^{1}A\tens_{A}E$ which obeys the condition
$\nabla(a.e)=da\tens e+a.\nabla e$
for all $e\in E$ and $a\in A$.
\end{defin}

If $E$ is an $A$-bimodule, we have the following idea of bimodule covariant derivative, see
\cite{DVMic,DVMass,MouLC,MadIntro}. 

\begin{defin} \label{ppll}  A {\em bimodule covariant derivative} on an
 $A$-bimodule $E$ is a triple $(E,\nabla,\sigma)$,
where $\nabla:E\to\Omega^{1}A\tens_{A}E$ 
is a left $A$-covariant derivative, and $\sigma:E\tens_A\Omega^1 A\to 
\Omega^1 A \tens_A E$ is a bimodule map obeying
\[
\nabla(e.a)\,=\,\nabla(e).a\,+\,\sigma(e\tens da)\ ,\qquad \forall\, e\in E,\ a\in A\ .
\]
\end{defin}

\section{Line modules}  \label{bhakvck}

To the reader who is familiar with Morita contexts (see Section~\ref{vcuadrts}), the material in this section will seem familiar. Propositions~\ref{cagjcjxdz}, \ref{bhcsavch} and \ref{kljhfvvhjg} are basically taken from  \cite{BassK}. 

\begin{defin}  \label{vzcvjhvkj}
Let $E$ be an $A$-bimodule that is finitely generated projective as a left $A$-module. 
If the coevaluation bimodule map $\coev:A\to E^\circ\tens_A E$ is an isomorphism,
then $E$ is called a \emph{weak left line module}. If in addition
the evaluation map $\ev:E\tens_A E^\circ\to A$ is an isomorphism,
then $E$ is called a \emph{left line} module.
\end{defin}

\begin{propos} \label{cagjcjxdz}
Let $E$ be an $A$-bimodule that is finitely generated projective as a left $A$-module. 
Then the following are equivalent:

\noindent $I$)\quad $E$ is a weak left line module.

\noindent $II$)\quad Every left module map from $E$ to $E$ is given by right action
by some element of $A$, and the only $a\in A$ for which $E.a=0$ is $a=0$.

\end{propos}
\noindent {\bf Proof:}\quad ($I\Rightarrow II$) First, suppose that
$E$ is a left line module and that $T:E\to E$ is a left module map. Then (summing over $i$)
$e_i\tens T(e^i)\in E^\circ\tens_A E$ and so is equal to $\coev(a)$ for some $a\in A$. Then
\begin{eqnarray*}
T(e) &=& T(e_i(e).e^i)\ =\ (\mathrm{ev}\tens\id)(e\tens \coev(a))\ =\ e.a\ .
\end{eqnarray*}
That $E.a=0$ implies $a=0$ follows from invertibility of the coevaluation.

($II\Rightarrow I$) Take $\beta\in E^\circ\tens_A E$. Then the map
\begin{eqnarray*}
e\mapsto (\mathrm{ev}\tens\id)(e\tens\beta)
\end{eqnarray*}
is a left $A$-module map from $E$ to $E$, and so is equal to $e.a$ for some fixed $a\in A$ by the hypothesis for ($II$). If we set $\gamma=\beta-\coev(a)\in E^\circ\tens_A E$, then
for all $e\in E$ we have $(\mathrm{ev}\tens\id)(e\tens\gamma)=0$. Now set
$\gamma=\sum_i \gamma_i\tens e^i$, and then we have for all $e\in E$,
(summed over $i$) $\gamma_i(e).e^i=0$. Applying $e_j$ to this gives 
$\gamma(e)\, e_j(e^i)=0$. Now
\begin{eqnarray*}
\gamma &=& \gamma_i\tens_A e^i \ =\ \gamma_i\tens_A e_j(e^i).e^j\cr
&=& \gamma_i.e_j(e^i) \tens_A e^j\ =\ 0\ .
\end{eqnarray*}
This shows that $\coev$ is surjective. To see that it is injective, its kernel is precisely those $a\in A$ for which $E.a=0$. \quad$\square$

\medskip There is a generalisation of Proposition~\ref{cagjcjxdz} part $II$ which will prove useful later. 

\begin{propos} \label{xdcxyjtalj}
Suppose that $E$ is a weak left line $A$-module.
Given a left $A$-module map $T:E\to E\tens_A F$ for some left $A$-module $F$, 
there is $f\in F$ so that $T(e)=e\tens f$, and $f$ is given by the image of $1\in A$ under the composition
\begin{eqnarray*}
A \stackrel{\coev} \longrightarrow E^\circ\tens_A E  \stackrel{\id\tens T} \longrightarrow 
E^\circ\tens_A E\tens_A F  \stackrel{\coev^{-1}\tens \id} \longrightarrow F\ .
\end{eqnarray*}
\end{propos}
\noindent {\bf Proof:}\quad 
From the definition of $f$,
\begin{eqnarray*}
(\coev\tens\id) (f) &=& (\id\tens T)\,(\coev(1))\ ,
\end{eqnarray*}
so
\begin{eqnarray*}
(\mathrm{ev}\tens\id)(e\tens(\coev\tens\id) f) &=& 
(\mathrm{ev}\tens\id)(e\tens(\id\tens T)\,(\coev(1)))\cr
&=& T((\mathrm{ev}\tens\id)(e\tens\coev(1))\cr
&=& T(e)\ .
\end{eqnarray*}
However we also have
\begin{eqnarray*}
(\mathrm{ev}\tens\id)(e\tens(\coev\tens\id) f) &=& e\tens f\ .\quad\square
\end{eqnarray*}

\medskip There is yet another version of this, which will be required!

\begin{propos}\label{hgcxaghjx}
Suppose that $E$ is a left line $A$-module.
Given a right $A$-module map $T:E\to F\tens_A E$ for some right $A$-module $F$, 
there is $f\in F$ so that $T(e)=f\tens e$, and $f$ is given by the image of $1\in A$ under the composition
\begin{eqnarray*}
A \stackrel{\ev^{-1}} \longrightarrow E\tens_A E^\circ  \stackrel{T\tens \id} \longrightarrow 
F \tens_A E\tens_A E^\circ  \stackrel{\id\tens\ev} \longrightarrow F\ .
\end{eqnarray*}
\end{propos}
\noindent {\bf Proof:}\quad By the definition of $f$,
\begin{eqnarray*}
(\id_F\tens \ev^{-1}\tens\id_E)(f\tens e) &=& (T\tens\id_{E^\circ}\tens\id_E)(\ev^{-1}\tens e)\ .
\end{eqnarray*}
From the equation $(\ev\tens\id_E)(\id_E\tens\coev)=\id_E$ we also have
\begin{eqnarray*}
(\id_F\tens \ev^{-1}\tens\id_E)(f\tens e) &=& (\id_F\tens\id_E\tens \coev)(f\tens e)\ .
\end{eqnarray*}
Now applying $\id_F\tens\id_E\tens\coev^{-1}$ to these equations gives the result.\quad
$\square$

\begin{propos} \label{bhcsavch}
 Suppose that $E$ is a weak left line $A$-module. Then
\begin{eqnarray*}
&&\mathrm{ev}\tens\id=\id\tens\coev^{-1}
:E\tens_A E^\circ\tens_A E\to E\ ,\cr
&&\id\tens\mathrm{ev}=\coev^{-1}\tens\id
:E^\circ\tens_A E\tens_A E^\circ\to E^\circ\ .
\end{eqnarray*}
\end{propos}
\noindent {\bf Proof:}\quad To show that the maps are equal in the first equation, we apply them after 
the invertible map
$\id\tens\coev:E\to E\tens_A E^\circ\tens_A E$. In the case of $\id\tens\coev^{-1}$, this just gives the identity. For the other case, we get
\begin{eqnarray*}
\big(\mathrm{ev}\tens\id\big)\big(\id\tens\coev\big):E\tens_A A\to A\tens_A E\ ,
\end{eqnarray*}
and this is the identity by the usual properties of evaluation and coevaluation.

To show that the maps are equal in the second equation, we apply them after 
the invertible map
$\coev\tens\id:E^\circ\to E^\circ\tens_A E\tens_A E^\circ$. 
In the case of $\coev^{-1}\tens\id$, this just gives the identity. For the other case, we get
\begin{eqnarray*}
\big(\id\tens\mathrm{ev}\big)\big(\coev\tens\id\big):E^\circ\tens_A A\to A\tens_A E^\circ\ ,
\end{eqnarray*}
and this is the identity by the usual properties of evaluation and coevaluation.
\quad$\square$

\begin{propos} \label{kljhfvvhjg}
 Suppose that $E$ is a weak left line $A$-module. If $\gamma$ is in the kernel of $\mathrm{ev}: E\tens_A E^\circ\to A$
 and $a\in \mathrm{image}(\mathrm{ev})$, then $a.\gamma=\gamma.a=0$. In particular, if
$\mathrm{ev}: E\tens_A E^\circ\to A$ is surjective, then it is an isomorphism. 
\end{propos}
\noindent {\bf Proof:}\quad Take $\gamma$ in the kernel of the evaluation map. Also take
$\beta\in E\tens_A E^\circ$ so that $\mathrm{ev}(\beta)=a$. Then by Proposition \ref{bhcsavch},
\begin{eqnarray*}
\gamma.a &=& (\id\tens\id\tens\mathrm{ev})(\gamma\tens\beta) \cr
&=& (\id\tens\coev^{-1}\tens\id)(\gamma\tens\beta) \cr
&=& (\mathrm{ev}\tens\id\tens\tens\id)(\gamma\tens\beta)\ =\ 0\ .
\end{eqnarray*}
The other way round is similar.\quad$\square$

\section{Morita contexts} \label{vcuadrts}
We begin with some definitions from \cite{BassK} (see also \cite{BerrKeat} for modern introduction to Morita theory). Suppose that $A$ and $B$ are unital algebras. 
Use ${}_A\mathbb{M}$ and ${}_B\mathbb{M}$ to denote the
category of left $A$- and left $B$-modules respectively, and ${}_A\mathbb{M}_B$ and ${}_B\mathbb{M}_A$  the $A-B$ and  $B-A$ bimodules
respectively. 

\begin{defin} \cite{BassK}\label{jkcvzhdshv}
A {\em Morita context} for $A$ and $B$ consists of $E\in {}_A\mathbb{M}_B$ and 
$F\in {}_B\mathbb{M}_A$ together with bimodule maps $\mu_1:E\tens_B F\to A$
and $\mu_2:F\tens_A E\to B$ so that
\begin{eqnarray*}
\mu_1\tens\id=\id\tens\mu_2&:&E\tens_B F\tens_A E \to E\ ,\cr
\mu_2\tens\id=\id\tens\mu_1&:&F\tens_A E\tens_B F \to F\ .
\end{eqnarray*}
This is called a {\em strict Morita context} in the case where $\mu_1$ and $\mu_2$ are surjective.
\end{defin}

It will be useful to recall the following results:

\begin{propos} \cite{BassK}\label{jkcvzhdshv1}
If $(A,B,E,F,\mu_1,\mu_2)$ is a strict Morita context as in Definition \ref{jkcvzhdshv}, then

\noindent a)\quad $\mu_1$ and $\mu_2$ are isomorphisms.

\noindent b)\quad $E$ and $F$ are finitely generated projective left $A$- and $B$-modules respectively.

\noindent c)\quad $E$ and $F$ are finitely generated projective right $B$- and $A$-modules respectively.
\end{propos}

The main reason for the interest in Morita contexts is given by the next result. Note that an equivalence of categories $\mathcal{C}$ and $\mathcal{D}$ is given by two functors
$P:\mathcal{C}\to\mathcal{D}$ and $Q:\mathcal{D}\to\mathcal{C}$ so that there are invertible natural 
transformations from $P\circ Q$ to the identity on $\mathcal{D}$, and from $Q\circ P$ to the identity on $\mathcal{C}$. 

\begin{propos} \cite{BassK}\label{jkcvzhdshv2}
There is a 1-1 correspondence between equivalences between the categories 
${}_A\mathbb{M}$ and ${}_B\mathbb{M}$ and strict Morita contexts
$(A,B,E,F,\mu_1,\mu_2)$. Corresponding to the strict Morita context, the functors are
given by
\begin{eqnarray*}
E\tens_B-& :& {}_B\mathbb{M}\to {}_A\mathbb{M}\ ,\cr
F\tens_A-& :& {}_A\mathbb{M}\to {}_B\mathbb{M}\ .
\end{eqnarray*}
\end{propos}

\section{Automorphisms of the centre} \label{cgfhjhfxsz}
Here we look at one of the implications of Proposition~\ref{cagjcjxdz}, which we will follow up later.
We use $Z(A)$ to denote the centre of the unital algebra $A$.  The reader should refer to
\cite{FroPicard}, where the map $\Phi$ and its algebraic implications are discussed.

\begin{propos} \label{cyjhgcxhj}
Given a weak left line module $L$, there is a unital algebra map $\Phi_L:Z(A)\to Z(A)$ given by
$z.e=e.\Phi_L(z)$ for all $a\in Z(A)$ and $e\in L$.
\end{propos}
\noindent {\bf Proof:}\quad If $z\in Z(A)$, then the map $e\mapsto z.e$ is a left $A$-module map from 
$L$ to $L$. By Proposition~\ref{cagjcjxdz} there is a $y\in A$ so that $z.e=e.y$ for all $e\in L$. 
By Proposition~\ref{cagjcjxdz} it follows that this $y$ is unique.
It follows that the map $e\mapsto e.y$ is a right $A$-module map, so
$e.a\,y=e.y\,a$ for all $e\in L$ and all $a\in A$. so $E.(a\,y-y\,a)=0$. By Proposition~\ref{cagjcjxdz} again, we see that 
$a\,y=y\,a$ for all $a\in A$, so $y\in Z(A)$. By uniqueness we get $\Phi_L(1_A)=1_A$. 

If $z,z'\in Z(A)$, then $z\,z'.e=z.e.\Phi_L(z')=e.\Phi_L(z)\,\Phi_L(z')$, so $\Phi_L(z\,z')=\Phi_L(z)\,\Phi_L(z')$. \quad$\square$

\medskip This map $\Phi_L$ is only interesting in the noncommutative case, as it measures the difference between the left and right module structures on $L$. 

\begin{propos}\label{cvuykxbcahisu}
If $L$ is a line module, then $\Phi_L:Z(A)\to Z(A)$ is invertible, with inverse
$\Phi_L^{-1}=\Phi_{L^\circ}$. 
\end{propos}
\noindent {\bf Proof:}\quad 
This uses the fact that $\coev:A\to L^\circ\tens_A L$ and $\ev:L\tens_A L^\circ\to A$ are bimodule maps. \quad$\square$

\begin{propos}\label{vhjcacvchj}
If $A$ is a star algebra and $L$ is a line module, then $\Phi_{\overline{L}}(z)=(\Phi_L^{-1}(z^*))^*$. 
\end{propos}
\noindent {\bf Proof:}\quad Compute
\begin{eqnarray*}
z.\overline{e}\ =\ \overline{e.z^*}\ =\ \overline{\Phi_L^{-1}(z^*).e}\ =\ \overline{e}.(\Phi_L^{-1}(z^*))^*\ .
\quad\square
\end{eqnarray*}

\begin{propos}  \label{cbvhjakjgcfgzx}
If weak line modules $L$ and $N$ are isomorphic as $A$-bimodules, then $\Phi_L=\Phi_N$. 
\end{propos}
\noindent {\bf Proof:}\quad Immediate.\quad$\square$

\section{The $\mathbb{N}$-graded and $\mathbb{Z}$-graded tensor algebras} \label{vcuskftsd}
For any $A$-bimodule $F$, we can define an $\mathbb{N}$-graded tensor algebra by
\begin{eqnarray}\label{jkhvjhvhkv}
T_{\mathbb{N}}(F)\ =\ A \oplus F \oplus (F\tens_A F) \oplus (F\tens_A F\tens_A F)\oplus\dots\ .
\end{eqnarray}
We write this as
\begin{eqnarray}
T_{\mathbb{N}}(F)^n\ =\ \left\{\begin{array}{cc}A & n=0 \\F^{\tens^n_A} & n>0\end{array}\right.\ .
\end{eqnarray}
The associative product on the algebra $T_{\mathbb{N}}(F)$ is just $\tens_A$. 
We shall also consider a $\mathbb{Z}$-graded object $T_{\mathbb{Z}}(F)$ by
\begin{eqnarray}\label{vhsvhbji}
T_{\mathbb{Z}}(F)^n\ =\ 
\left\{\begin{array}{cc}A & n=0 \\F^{\tens^n_A} & n>0 \\(F^\circ)^{\tens^{-n}_A} & n<0\end{array}\right.\ .
\end{eqnarray}
In general this will not form an algebra under $\tens_A$ as $F\tens_A F^\circ\neq A$. However we can take a special case:

\begin{propos}\label{jkhcvhjcjhcg}
If $L$ is a weak left line $A$-bimodule, then $T_{\mathbb{Z}}(L)$ with product $\tens_A$ combined with $\ev$ and $\coev^{-1}$
gives an associative algebra. 
\end{propos}

Before we prove Proposition \ref{jkhcvhjcjhcg}, we should comment that in the $C^*$-algebra case
(with account taken of norms), it is likely that this construction gives the Pimsner algebra associated to the bimodule $L$
(see \cite{PimsAlg,DopPinZuc,PinRob}), as the Pimsner algebra is constructed as a representation on a space which is just $T_{\mathbb{N}}(F)$ with `creation operators' acting by $\tens$. 

We need to be more precise about the product in Proposition~\ref{jkhcvhjcjhcg}. 
For notation we set $L^0=A$ and
 for $n\ge 1$ we set $L^n=L^{\tens^n_A}$ and $L^{-n}=(L^\circ)^{\tens^{-n}_A}$. Now recursively define, for $n\ge 1$,
 \begin{eqnarray}\label{bvhsvbhk1}
\ev^n:L^n\tens_A L^{-n}\to A\ ,
\end{eqnarray}
by $\ev^1=\ev$ and $\ev^{n+1}=\ev\, (\id\tens\ev^n\tens\id)$. Also define, for $n\ge 1$,
 \begin{eqnarray}\label{bvhsvbhk2}
\tilde\ev^n:L^{-n}\tens_A L^{n}\to A\ ,
\end{eqnarray}
by $\tilde\ev^1=\coev^{-1}$ and $\tilde\ev^{n+1}=\coev^{-1}\, (\id\tens\tilde\ev^n\tens\id)$. 

The product $L^n\tens_A L^m\to L^{n+m}$ is simply $\tens_A$ in the case where
$n$ and $m$ have the same sign. For different signs, we use the isomorphisms
$\ev^n$ and $\tilde\ev^n$. We consider the possible cases, for $n,m>0$:
\begin{eqnarray*}
\tilde\ev^n\tens\id:L^{-n}\tens_A L^m\to L^{m-n} & \mathrm{for} & m\ge n\ ,\cr
\id\tens\tilde\ev^m:L^{-n}\tens_A L^m\to L^{m-n} & \mathrm{for} & m\le n\ ,\cr
\ev^n\tens\id:L^{n}\tens_A L^{-m}\to L^{n-m} & \mathrm{for} & m\ge n\ ,\cr
\id\tens\ev^m:L^{n}\tens_A L^{-m}\to L^{n-m} & \mathrm{for} & m\le n\ .
\end{eqnarray*}

\smallskip
\noindent {\bf Proof of Proposition \ref{jkhcvhjcjhcg}:}\quad We need to compare the products on
$(L^n\tens_A L^m)\tens_A L^r$ and  $L^n\tens_A (L^m\tens_A L^r)$. By associativity of
$\tens_A$, these give the same result if $n,m,r$ have the same sign. In fact it is not very difficult to see that we obtain exactly the same result if $n,m$ have the same sign, and if $m,r$ have the same sign. This leaves the more difficult cases where the sequence $n,m,r$ alternates in sign. The case where the alternating sequence $n,m,r$ has numbers $\pm 1$ is exactly Proposition \ref{bhcsavch}. 
The problem is to extend this to other alternating sequences. We will suppose that $n,r>0$ and $m<0$ (we write $m=-s$ for convenience) -- the other case is almost identical. For clarity,
we give $\id^q$ for the identity on either $L^q$ or $L^{-q}$. 
The first mentioned product on equations (\ref{jkjkjkjk1}--\ref{jkjkjkjk5}) is that on 
$(L^n\tens_A L^{-s})\tens_A L^r$, and the second is that on $L^n\tens_A (L^{-s}\tens_A L^r)$. 

\noindent {Case 1:}\quad $s\ge n+r$. In this case the required equality between products is
\begin{eqnarray} \label{jkjkjkjk1}
\ev^{n}\tens \id^{s-n-r}\tens\tilde\ev^r\ =\ \ev^{n}\tens \id^{s-n-r}\tens\tilde\ev^r :
L^n\tens_A L^{-s}\tens_A L^r\to L^{n+r-s}\ .
\end{eqnarray}

\noindent {Case 2:}\quad $n,r\ge s$. In this case the required equality between products is
\begin{eqnarray}\label{jkjkjkjk2}
\id^{n-s}\tens \ev^s\tens\id^r\ =\ \id^n\tens\tilde\ev^s\tens\id^{r-s}:
L^n\tens_A L^{-s}\tens_A L^r\to L^{n+r-s}\ .
\end{eqnarray}
By removing identities from both sides, it is enough to prove this for $n=r=s$. 

\noindent {Case 3:}\quad $n\ge s\ge r$. In this case the required equality between products is
\begin{eqnarray}\label{jkjkjkjk3}
\id^{n-s}\tens \ev^s\tens\id^r\ =\ \id^{n-s+r}\tens\ev^{s-r}\tens\tilde\ev^{r}:
L^n\tens_A L^{-s}\tens_A L^r\to L^{n+r-s}\ .
\end{eqnarray}
By removing identities from the left, it is enough to prove this for $n=s\ge r$, i.e.\ that
\begin{eqnarray*}
\ev^n\tens\id^r\ =\ \id^{r}\tens\ev^{n-r}\tens\tilde\ev^{r}:
L^n\tens_A L^{-n}\tens_A L^r\to L^{r}\ .
\end{eqnarray*}
But both of these sides can be given as a composition with the same beginning, as
\begin{eqnarray*}
\ev^n\tens\id^r &=& (\ev^r\tens\id^r)\, (\id^{r}\tens\ev^{n-r}\tens\id^{2\,r})\ ,\cr
\id^{r}\tens\ev^{n-r}\tens\tilde\ev^{r} &=& (\id^r\tens \tilde\ev^{r})\, (\id^{r}\tens\ev^{n-r}\tens\id^{2\,r})\ .
\end{eqnarray*}
Using this we are reduced to the case $n=r=s$ again.

\noindent {Case 4:}\quad $n\ge s\ge r$. In this case the required equality between products is
\begin{eqnarray}\label{jkjkjkjk4}
\ev^{n}\tens \tilde\ev^{s-n}\tens\id^{n+r-s}\ =\ \id^{n}\tens\tilde\ev^{s}\tens\id^{r-s}:
L^n\tens_A L^{-s}\tens_A L^r\to L^{n+r-s}\ .
\end{eqnarray}
By removing identities from the right, it is enough to prove this for $r=s\ge r$, i.e.\ that
\begin{eqnarray*}
\ev^{n}\tens \tilde\ev^{s-n}\tens\id^{n}\ =\ \id^{n}\tens\tilde\ev^{s}:
L^n\tens_A L^{-n}\tens_A L^r\to L^{r}\ .
\end{eqnarray*}
But both of these sides can be given as a composition with the same beginning, as
\begin{eqnarray*}
\ev^n\tens\id^r &=& (\ev^n\tens\id^n)\, (\id^{2\,n}\tens\tilde\ev^{s-n}\tens\id^{n})\ ,\cr
\id^{r}\tens\ev^{n-r}\tens\tilde\ev^{r} &=& (\id^n\tens \tilde\ev^{n})\, (\id^{2\,n}\tens\tilde\ev^{s-n}\tens\id^{n})\ .
\end{eqnarray*}
Using this we are reduced to the case $n=r=s$ again.

\noindent {Case 5:}\quad $n+r\ge s\ge n,r$. In this case the required equality between products is
\begin{eqnarray}\label{jkjkjkjk5}
\ev^{n}\tens \tilde\ev^{s-n}\tens\id^{r+n-s}\ =\ \id^{n-s+r}\tens\ev^{s-r}\tens\tilde\ev^{r}:
L^n\tens_A L^{-s}\tens_A L^r\to L^{n+r-s}\ .
\end{eqnarray}
But both of these sides can be given as a composition with the same beginning, as
\begin{eqnarray*}
\ev^{n}\tens \tilde\ev^{s-n}\tens\id^{r+n-s} &=& (\ev^{n+r-s}\tens\id^{n+r-s})
\cr &&\, (\id^{n+r-s}\tens\ev^{s-r}\tens\id^{r+n-s}\tens\tilde\ev^{s-n}\tens\id^{r+n-s})\ ,\cr
\id^{n-s+r}\tens\ev^{s-r}\tens\tilde\ev^{r} &=& (\id^{n+r-s}\tens\tilde\ev^{n+r-s})
\cr &&\, (\id^{n+r-s}\tens\ev^{s-r}\tens\id^{r+n-s}\tens\tilde\ev^{s-n}\tens\id^{r+n-s})\ . 
\end{eqnarray*}
Using this we are reduced to the case $n=r=s$ again.

\smallskip All but one of these cases reduces to showing the following equality
\begin{eqnarray} \label{jkjkjkjk10}
\ev^{n}\tens \id^{n}\ =\ \id^{n}\tens\tilde\ev^n :
L^n\tens_A L^{-n}\tens_A L^n\to L^{n}\ .
\end{eqnarray}
We do this by induction on $n$, noting that the starting case $n=1$ is done in Proposition 
\ref{bhcsavch}. Now suppose that (\ref{jkjkjkjk10}) is true. We need to show the following 
equality:
\begin{eqnarray} \label{jkjkjkjk11}
\ev^{n+1}\tens \id^{n+1}\ =\ \id^{n+1}\tens\tilde\ev^{n+1} :
L^{n+1}\tens_A L^{-n-1}\tens_A L^{n+1}\to L^{n+1}\ .
\end{eqnarray}
Then, by using Proposition \ref{bhcsavch},
\begin{eqnarray*}
\ev^{n+1}\tens \id^{n+1} &=& (\ev\tens \id^{n+1})\, (\id\tens\ev^{n}\tens \id^{n+2}) \cr
&=& (\id\tens\tilde\ev\tens \id^{n})\, (\id\tens\ev^{n}\tens \id^{n+2})\cr
&=& (\id\tens\ev^{n}\tens \id^{n})\, (\id^{2\,n+1}\tens\tilde\ev\tens\id^{n})\ ,
\end{eqnarray*}
and using (\ref{jkjkjkjk10}) gives
\begin{eqnarray*}
\ev^{n+1}\tens \id^{n+1} &=& 
 (\id^{n+1}\tens\tilde\ev^{n})\, (\id^{2\,n+1}\tens\tilde\ev\tens\id^{n})\cr
 &=& \id^{n+1}\tens\tilde\ev^{n+1}\ .\quad\square
\end{eqnarray*}

\section{Hopf-Galois extensions} \label{cytjadsdrta}
The reader is reminded that if $H$ is a Hopf algebra and $B$ a unital algebra, then
$B$ is a comodule algebra for the right action $\rho:B\to B\tens H$
(written $\rho(b)=b_{[0]}\tens b_{[1]}$) if the following
equations hold for all $b,b'\in B$:
\begin{eqnarray*}
\rho(b\,b')\ =\ b_{[0]}\, b'_{[0]}\tens b_{[1]}\,b'_{[1]}\ ,\quad \rho(1_B)\ =\ 1_B\tens 1_H\ .
\end{eqnarray*}
If $A$ is the invariant subalgebra of $B$ (i.e.\ $\{b\in B|\rho(b)=b\tens 1_H\}$),
 then the canonical map is
\begin{eqnarray*}
\mathrm{can}:B\tens_A B\to B\tens H\ ,\quad b'\tens b\mapsto b'\,b_{[0]}\tens b_{[1]}\ .
\end{eqnarray*}
If the canonical map is a 1-1 correspondence, then $B$ is said to be a 
\emph{Hopf-Galois extension} of $A$.

Any $\mathbb{Z}$-graded algebra $C=\oplus_{n\in \mathbb{Z}} C_n$ can be thought of as a comodule algebra over $\mathbb{CZ}$, the group algebra of $(\mathbb{Z},+)$.
This is given by mapping $c\in C_n$ to $ c\tens\underline{n}$, where we write $\underline{n}$
to be a generator of the group algebra of $\mathbb{Z}$. The invariant subalgebra is then $C_0$. 

\begin{propos}  \label{kscvkcvvh}
The  $\mathbb{Z}$-graded algebra $C=\oplus_{n\in \mathbb{Z}} C_n$ is a 
$\mathbb{CZ}$ Hopf-Galois extension of $C_0$ if and only if every product $:C_n\tens C_m\to C_{n+m}$ is surjective. 
\end{propos}
\noindent {\bf Proof:}\quad For the $\mathbb{CZ}$ coaction, the canonical map
can be written as
\begin{eqnarray*}
c\tens c' \mapsto c\,c'\tens\underline{m}
\end{eqnarray*}
for $c'\in C_m$. For the canonical map to be surjective, we must have the map
$C\tens C_m\to C$ being surjective, so it follows that $:C_n\tens C_m\to C_{n+m}$ is surjective
for all $n,m\in\mathbb{Z}$. \quad$\square$

\medskip In fact Proposition~\ref{kscvkcvvh} holds more generally for any algebra which is graded by a group, as is shown in Prop.\ 8.2.1 of \cite{caenBook}. 
 
 \begin{propos}\label{cvbhcvv}
 For a weak left line module $L$ over the unital algebra $A$, the tensor algebra
 $T_{\mathbb{Z}}(L)$ is a $\mathbb{CZ}$ Hopf-Galois extension of $A$ if and only if $L$ is a
 left line module. 
\end{propos}
\noindent {\bf Proof:}\quad We need to show that every product $:T_{\mathbb{Z}}(F)_n\tens T_{\mathbb{Z}}(F)_m\to T_{\mathbb{Z}}(F)_{n+m}$ is surjective. We split this into cases: First note that it is automatic for either of $n$ or $m$ being zero. If $n,m>0$ we have 
\begin{eqnarray*}
\tens_A: L^{\tens^n_A} \tens L^{\tens^m_A} \to L^{\tens^{n+m}_A} \ ,
\end{eqnarray*}
which is onto by definition. Likewise if $n,m<0$ we have the surjective map
\begin{eqnarray*}
\tens_A: (L^\circ)^{\tens^{-n}_A} \tens (L^\circ)^{\tens^{-m}_A} \to (L^\circ)^{\tens^{-n-m}_A}\ .
\end{eqnarray*}
It is not too hard to see that the remaining cases, where $n$ and $m$ have opposite sign,
reduce to the cases of showing that the following two maps are onto, for all $r>0$:
\begin{eqnarray*}
\ev^r:L^r\tens_A L^{-r}\to A\ \mathrm{and}\ \tilde\ev^r:L^{-r}\tens_A L^{r}\to A\ .
\end{eqnarray*}
By the recursive definitions of $\ev^r$ and $\tilde\ev^r$ this follows immediately from the case $r=1$, 
and this is the case if $L$ is a left line module. \quad$\square$

\begin{theorem} \label{vchuaytdrtys}
Let $A$ be an algebra.
There is a 1-1 correspondence between:

\noindent a)\quad Autoequivalences of the category ${}_A\mathbb{M}$ of left $A$-modules.

\noindent b)\quad Left line modules over $A$.

\noindent c)\quad Hopf-Galois $\mathbb{CZ}$ (group algebra of $\mathbb{Z}$) extensions of $A$. 
\end{theorem}
\noindent {\bf Proof:}\quad ($a\Rightarrow b$)\quad From Proposition \ref{jkcvzhdshv2}, autoequivalences of the category ${}_A\mathbb{M}$
correspond to strict Morita contexts $(A,A,E,F,\mu_1,\mu_2)$. 
By Proposition \ref{jkcvzhdshv1}, $E$ is a finitely generated projective left $A$-module, and
we have $A$-bimodule isomorphisms $\mu_1:E\tens_A F\to A$
and $\mu_2:F\tens_A E\to A$ obeying the conditions of \ref{jkcvzhdshv}. If we relabel the maps
as $\ev=\mu_1:E\tens_A F \to A$ and $\coev=\mu_2^{-1}:A\to F\tens_A E$, we see that $F=E^\circ={}_A\End(E,A)$. 
The `associativity' conditions in \ref{jkcvzhdshv} ensure the usual behaviour of the evaluation and coevaluation maps. Then $E$ satisfies the conditions to be a left line module given in \ref{vzcvjhvkj}.

($b\Rightarrow c$)\quad Given a left line module $L$, $T_{\mathbb{Z}}(L)$ is a 
$\mathbb{CZ}$ Hopf-Galois extension of $A$ by Proposition \ref{cvbhcvv}.

($c\Rightarrow a$)\quad Given an algebra $C$ which is a $\mathbb{CZ}$ Hopf-Galois extension of $A$, the $\mathbb{CZ}$ coaction splits $C$ into a direct sum of integer graded parts $C_n$. 
We have $C_0=A$, and set $E=C_1$ and $F=C_{-1}$. Then $E$ and $F$ are $A$-bimodules,
and the multiplication maps $\mu_1:E\tens_A F\to A$ and $\mu_2:F\tens_A E\to A$
are onto by the Hopf-Galois condition in Proposition \ref{kscvkcvvh}. The `associativity' conditions in
\ref{jkcvzhdshv} are implied by the associativity of the algebra $C$.\quad$\square$

\section{The Picard group of an algebra} \label{pichgcjf}
In topology the Picard group $\mathrm{Pic}(X)$ has elements the line bundles on a topological space $X$ (up to isomorphism), and the group product is tensor product. In the noncommutative case, suppose that we have two line modules, $L$ and $M$, over an algebra $A$. Then $L\tens_A M$ is another line module, as we now show. Define $(L\tens_A M)^\circ=M^\circ\tens_A L^\circ$, with 
\begin{eqnarray*}
\ev_{L\tens M} &=&  \ev_L(\id_L\tens\ev_M\tens\id_{L^\circ}  ) : L\tens_A M\tens_A M^\circ\tens_A L^\circ \to A\ ,\cr
\coev_{L\tens M} &=&  (\id_{M^\circ}\tens\coev_L\tens\id_{M}  ) \coev_M : A\to M^\circ\tens_A L^\circ\tens_A L\tens_A M\ .
\end{eqnarray*}
We could define $\mathrm{Pic}(A)$ for an algebra $A$ to be the isomorphism classes of line modules under $\tens_A$ (see \cite{BassK,GAMS,Bolla,Masuoka1}).

Now we give a result on the algebra map of the centre given in Proposition~\ref{cyjhgcxhj}

\begin{propos} \label{cyjhgcxhjhgi} 
The map $\Phi:\mathrm{Pic}(A)\to \mathrm{Aut}(Z(A))$ given by the isomorphism class of $L$ mapping to $\Phi_L$ is an order reversing group homomorphism. 
\end{propos}
\noindent {\bf Proof:}\quad For $e\tens f\in L\tens_A N$ we have
\begin{eqnarray*}
z.e\tens f \ =\ e.\Phi_L(z)\tens f \ =\ e\tens \Phi_L(z).f \ =\ e\tens f.\Phi_N(\Phi_L(z)) \ .
\quad\square
\end{eqnarray*}

Suppose we wished to look at this from the point of view of principal bundles, or rather $\mathbb{Z}$
Hopf-Galois extensions. We might ask if $T_{\mathbb{Z}}(L\tens_A M)$ was a subalgebra of
$T_{\mathbb{Z}}(L) ?? T_{\mathbb{Z}}(M)$, where $??$ denotes some sort of product. We know of no answer to this, but we have some comments. One way to get such a product would be to suggest some form of `exchange law' (e.g.\ see results on factorisation in Section 7 of \cite{Ma:Book} or \cite{CapSchVan}), possibly a bimodule map from $L\tens_A M$ to $M\tens_A L$. If such a bimodule map existed and was invertible, then $L$ and $M$ would commute in 
$\mathrm{Pic}(A)$. We do not know if $\mathrm{Pic}(A)$ is abelian in general (though we would not be surprised if someone did know). One way to get such bimodule maps would be to look at Section 4.2 of 
\cite{BegMa2}, where a star operation on all of $L$, $M$ and $L\tens_A M$ is used to give
an invertible map $\Gamma:L\tens_A M\to M\tens_A L$. Of course another comment on this is that it is not at all obvious how to extend a star operation to tensor products. 

Following this, it is tempting to think that there is some sort of universal extension of $A$ for line bundles, which would be a $\mathrm{Pic}(A)$-graded Hopf-Galois extension of $A$, for which the $[L]$  degree part of the algebra is isomorphic to $L$ as a bimodule. The problem here is that it is not at all obvious that an associative product can be derived if we have to choose representatives of isomorphism classes of bimodules at every stage. This is then another open problem. Since the quotient of two bimodule isomorphisms between line modules is an invertible element of the centre $Z(A)$ (from Proposition~\ref{cagjcjxdz}), it is quite likely that cocycles will be involved.

\section{Hermitian metrics on modules} \label{gahcvjdvyu}
A Hermitian metric is the equivalant of a Riemannian metric, but on an arbitrary complex vector bundle. It is similar to the idea of Hilbert $C^*$-module in $C^*$-algebra theory (see \cite{Lance}).
As is usual on complex vector spaces, the (generalised) inner product is linear in one variable, and conjugate linear in the other. At the moment, we do not make any assumption about positivity.

\begin{defin}\label{hermdeff}\cite{BegMa4} 
A \emph{non degenerate Hermitian structure} on an $A$-bimodule $E$ is given by an invertible
$A$-bimodule morphism
$G:\overline{E}\to E^\circ$. From this we define an \emph{inner product}
$\<,\>=\ev_E(\id\tens G):E\tens_A \overline{E}\to A$, and this is required to satisfy the condition
that the following composition is just $\<,\>$:
\begin{eqnarray*}
E \tens_A \overline{E} \stackrel{\mathrm{bb}\tens\id} \longrightarrow
\overline{ \overline{E}} \tens_A \overline{E} \stackrel{\Upsilon^{-1}} \longrightarrow
\overline{E \tens_A \overline{E}} \stackrel{\overline{\<,\>}} \longrightarrow \bar A 
\stackrel{\star^{-1}} \longrightarrow  A \ .
\end{eqnarray*}
\end{defin}

We write $\<e,\bar f\>=\ev(e \tens_A G(\bar f))$. 
For $e,f\in E$ the composition in \ref{hermdeff}  is
\begin{eqnarray*}
e\tens \bar f \longmapsto  \overline{\bar e} \tens \bar f 
\longmapsto \overline{f\tens\bar e} \longmapsto 
 \overline{\< f,\bar e\>} \longmapsto \< f,\bar e\>^*\ .
\end{eqnarray*}
Thus the condition in \ref{hermdeff} is that $\<e,\bar f\>=\<f,\bar e\>^*$. There are some other
formulae which are virtually automatic from the definition. 
Since $\ev_E$ is a right $A$-module map, for all $a\in A$,
\begin{eqnarray*}
\<e,\overline{a.f}\>\,=\, \<e,\bar f.a^*\>\,=\, \<e,\bar f\>\,a^*\ .
\end{eqnarray*}
Since $\ev_E$ is a left $A$-module map, for all $a\in A$,
\begin{eqnarray*}
\<a.e,\bar f\>\,=\, a\,\<e,\bar f\>\ .
\end{eqnarray*}
Since we are using the tensor product over $A$,
\begin{eqnarray*}
\<e.a,\bar f\>\,=\, \<e,a.\bar f\>\,=\, \<e,\overline{f.a^*}\>\ .
\end{eqnarray*}

\medskip The following proposition from \cite{BegMa4} will prove useful.

\begin{propos}  \label{xsrhkldjtyfggh}
Suppose that $E$ is finitely generated projective as a left module, with dual basis
$e_i\tens e^i\in E^\circ\tens E$, and 
let $G$ be a non-degenerate Hermitian structure on $E$. 
Set $g^{ij}=\<e^i,\overline{e^j}\>$, so it is automatic that
$g^{ij*}=g^{ji}$. Then  $G(\overline{e^i})=e_j.g^{ji}$ 
(summation convention applies). Define 
$G^{-1}(e_i)=\overline{g_{ij}.e^j}$, where without
loss of generality
it can be assumed that $g_{ij}.\ev(e^j\tens e_k)= g_{ik}$. Then:

a)\quad $g^{ij}\,g_{jk}\,=\, \ev(e^i\tens e_k)$\ .

b)\quad $g_{ij}\,g^{jk}\,=\, \ev(e^k \tens e_i)^*$\ .

c)\quad $g_{iq}^*\,=\, g_{qi}$\ .

\end{propos}\label{cvgads}
\proof Begin with
\begin{eqnarray*}
\overline{e^i} \,=\, G^{-1}(G(\overline{e^i}))\,=\, G^{-1}(e_j.g^{ji}) \,=\, G^{-1}(e_j) .g^{ji}\,=\, 
\overline{g_{jk}.e^k} .g^{ji}  \,=\, 
\overline{g^{ij}\,g_{jk}.e^k}    \  ,
\end{eqnarray*}
and apply $e_n$ to both sides to get (a). Also
\begin{eqnarray*}
e_i \,=\, G(G^{-1}(e_i))\,=\, G(\overline{g_{ij}.e^j})\,=\, 
G(\overline{e^j}).g_{ij}^*\,=\, e_k.g^{kj}\, g_{ij}^*\ ,
\end{eqnarray*}
and applying both sides to $e^p$ gives $\ev(e^k \tens e_i)=g^{kj}\, g_{ij}^*$,
while applying $*$ to this gives
(b). Finally
\begin{eqnarray*}
g_{ni}\,=\,
g_{nk}\,\ev(e^k \tens e_i) \,=\, g_{nk}\,g^{kj}\, g_{ij}^*\,=\,\ev(e^j \tens e_n)^*\, g_{ij}^*
\,=\,(g_{ij}\, \ev(e^j \tens e_n))^*\,=\,g_{in}^*\ .\quad\square
\end{eqnarray*}

\medskip Note that the fact that we have not defined $g_{nk}$ as the inverse to the matrix
$g^{nk}$ is nothing to do with noncommutativity. Even in ordinary differential geometry, this identification with the inverse requires choosing a chart which trivialises the bundle.

In the presence of a nondegenerate Hermitian metric, we have the following corollary to the results of Section~\ref{cgfhjhfxsz}:

\begin{cor}
If $A$ is a star algebra and $L$ is a line module with nondegenerate Hermitian metric, then
$\Phi_L:Z(A)\to Z(A)$ is a star algebra map. 
\end{cor}
\noindent {\bf Proof:}\quad A nondegenerate Hermitian metric on $L$ is simply an
$A$-bimodule isomorphism between $\overline{L}$ and $L^\circ$. Now use  
Propositions \ref{cvuykxbcahisu}, \ref{vhjcacvchj} and \ref{cbvhjakjgcfgzx}.\quad$\square$

\section{Involution on the $\mathbb{Z}$-graded tensor algebra} \label{vuwicvhcv}
For the $\mathbb{Z}$-graded tensor algebra we will give a star operation
which sends (the conjugate of) grade $1$ to  grade $-1$, i.e.\ 
$\overline{L}$ and $L^\circ$ are related by the star map.  
There is an example from complex geometry to motivate
this star operation:
The canonical line bundle in algebraic geometry has a star operation, given by conjugation of the complex coordinates, but that operation takes values in the dual of the canonical bundle rather than the original bundle. 
The reader should note that specifying a bimodule isomorphism 
between $\overline{L}$ and $L^\circ$ is just giving a non-degenerate Hermitian structure to the module $L$.

Given a non-degenerate Hermitian inner product $G: \overline{L}\to L^\circ$ on $L$, we define a star operation $z\mapsto
\overline{z^*}$ by
\begin{eqnarray}\label{jhvshjvhvch}
L\stackrel{\mathrm{bb}}\longrightarrow \overline{ \overline{L}}
\stackrel{\overline{G}}\longrightarrow \overline{ L^\circ}\ ,\quad
L^\circ \stackrel{G^{-1}}\longrightarrow \overline{L}\ .
\end{eqnarray}
In terms of the dual basis and the matrices in Proposition~\ref{xsrhkldjtyfggh} this means
\begin{eqnarray}\label{uvvvgcrxe}
e^{i*}\ =\ e_j\,g^{ji}\ ,\quad {e_k}^*\ =\ g_{kq}\, e^q\ .
\end{eqnarray}
Now we have to check the compatibility of the star operation and the product in 
$T_{\mathbb{Z}}(L)$:

\begin{lemma}  \label{xfghkjhxd}
For all $x\in L$ and $y\in L^\circ$ in $T_{\mathbb{Z}}(L)$, we have, using the star operation 
in (\ref{jhvshjvhvch}),
$(x\,y)^*=y^*\,x^*$ and $(y\,x)^*=x^*\,y^*$. 
\end{lemma}
\noindent {\bf Proof:}\quad First we verify $(x\,y)^*=y^*\,x^*$. The two maps from
$L\tens_A L^\circ$ to $A$ given by $x\tens y\mapsto x\,y$ and 
$x\tens y\mapsto (y^*\,x^*)^*$ are bimodule maps, so it is sufficient to verify that
they are equal on a left basis of $L$ and a right basis of $L^\circ$, so we use
$x=e^i$ and $y=e_k$. Then, by (\ref{uvvvgcrxe}), where $P_{qj}=\ev(e^q\tens e_j)$,
\begin{eqnarray*}
y^*\,x^* &=& g_{kq}\, e^q\, e_j\,g^{ji}\ =\ g_{kq}\, P_{qj}\,g^{ji}\cr
&=& g_{kq}\, g^{qi}\ =\ (P_{ik})^*\ =\ (e^i\, e_k)^*\ .
\end{eqnarray*}
Verifying $(y\,x)^*=x^*\,y^*$ is rather more difficult. If we set $a=x^*\,y^*$, by definition
\begin{eqnarray*}
\coev(a)\ =\ x^*\tens y^*\in L^\circ\tens_A L\ .
\end{eqnarray*}
This means that, for every basis element $e^n\in L$,
\begin{eqnarray*}
e^n.a &=& \ev(e^n\tens x^*).y^*\ .
\end{eqnarray*}
Applying star to this gives
\begin{eqnarray*}
a^*\,e_m\,g^{mn} &=& y.\ev(e^n\tens x^*)^*\ .
\end{eqnarray*}
Using the formulae for the matrices, we get (using the previous result for the last step)
\begin{eqnarray*}
a^*\,e_m &=& y.\ev(e^n\tens x^*)^*\,g_{nm}\cr
&=& y.\ev(g_{mn}\,e^n\tens x^*)^* \cr
&=& y.\ev({e_m}^*\tens x^*)^*\cr
&=& y.\ev(x\tens e_m)\ .
\end{eqnarray*}
Then
\begin{eqnarray*}
\ev(e^p\tens a^*\,e_m) &=& \ev(e^p\tens y)\, \ev(x\tens e_m)\ ,
\end{eqnarray*}
which we rearrange to obtain
\begin{eqnarray*}
\ev(e^p\,a^*\tens e_m) &=& \ev\big(\ev(e^p\tens y)\,x\tens e_m\big)\ .
\end{eqnarray*}
As this is true for all $e^m$, we get, for all $e^p$,
\begin{eqnarray*}
e^p\,a^* &=& \ev(e^p\tens y)\,x\ .
\end{eqnarray*}
But this is exactly the statement that $\coev(a^*)=y\tens x\in L^\circ\tens_A L$.
\quad$\square$

\begin{propos}
If  the star operation on $T_{\mathbb{Z}}(L)$ is defined in terms of the previous star operations on $L$ and $L^\circ$ by order reversal,
\begin{eqnarray*}
x_1\tens x_2\tens\dots\tens x_m &\longmapsto & x_m^*\tens \dots\tens x_2^*\dots\tens x_1^*\ ,
\end{eqnarray*}
then $T_{\mathbb{Z}}(L)$ becomes a star algebra. 
\end{propos}
\noindent {\bf Proof:~} There are only two non-trivial products to check consistency with,
that is $L\tens_A L^\circ\to A$ and $L^\circ\tens_A L\to A$. This is done in Lemma~\ref{xfghkjhxd}.
\quad$\square$

\section{Involution on the $\mathbb{N}$-graded tensor algebra} \label{xctkvdft}
Here we assume that there is a star operation on a bimodule $E$, which we write
as a bimodule map
$\star:E\to \overline{E}$ or $e\mapsto \overline{e^*}$. We stress that this is an entirely different construction to that for the $\mathbb{Z}$-graded tensor algebra, we are not assuming the existence of a Hermitian metric on $E$. 

Given this, we define a star operation $\star_n$ on $E^{\tens_A^m}$ by
\begin{eqnarray*}  \label{vcsvcksv}
x_1\tens x_2\tens\dots\tens x_m &\longmapsto & x_m^*\tens \dots\tens x_2^*\dots\tens x_1^*\ .
\end{eqnarray*}
The algebra $T_{\mathbb{N}}(E)$ with this star operation becomes a star algebra. Unlike the $T_{\mathbb{Z}}(E)$ case, there are no non-trivial relations to check. 

 In general it is not expected that $T_{\mathbb{N}}(E)$ will be the star algebra of functions on the total space of a bundle. It is too big, as we have not included any form of symmetry requirement. [Recall that on $\mathbb{R}^2$, the functions are polynomials in the coordinates, and therefore obey a symmetry condition.] However in the case of a line bundle, this is not a concern, as classically locally there
is only one generator of the functions on the bundle - the fibre coordinate. If $L$ is a line module,
we should think of $T_{\mathbb{N}}(L)$ as the star algebra of functions on the line bundle which are polynomial in the fibre direction.

\section{The $\mathbb{Z}/2$-graded Hopf-Galois extension} \label{cvaufv2}

In the $\mathbb{N}$-graded tensor algebra $T_{\mathbb{N}}(L)$, the
summand $L^{\tens_A^n}$ can be thought of as the functions on the line bundle
with polynomial growth of order $n$. 
In this section we suppose that $L$ has a non-degenerate Hermitian metric $G:\overline{L}\to L^\circ$ and a star operation $\star:L\to\overline{L}$. 

\begin{defin}\label{cfuadikf}
Define $\xi\in L\tens_A L$ to be the image of $1_A$ under
\begin{eqnarray*}
A\stackrel{\coev}  \longrightarrow L^\circ\tens_A L
 \stackrel{G^{-1}\tens\id}   \longrightarrow   \overline{L}\tens_A L
  \stackrel{\star^{-1}\tens\id}   \longrightarrow   L\tens_A L\ .
\end{eqnarray*}
Also define $\alpha\in A$ to be the image of $1_A$ under
\begin{eqnarray*}
A \stackrel{\coev} \longrightarrow L^\circ\tens_A L 
\stackrel{(G\,\star)^{-1}\tens G\,\star} \longrightarrow L\tens_A L^\circ 
\stackrel{\ev} \longrightarrow A\ .
\end{eqnarray*}
\end{defin}

\begin{propos}  \label{cfgyajkhj}
The element $\alpha\in A$ in Definition~\ref{cfuadikf} is invertible and central in $A$.
 Similarly $a.\xi=\xi.a$ for all $a\in A$.  Also $\xi$ and $\alpha$ are Hermitian,
 as can be deduced from the following formulae:
\begin{eqnarray*}
\alpha &=& \<e^{j*}\,g_{ji},\overline{e^{i*}}\>\ ,\cr
\xi &=& e^{j*}\,g_{ji}\tens e^i\ .
\end{eqnarray*}
\end{propos}
\noindent {\bf Proof:}\quad As $a.1_A=a=1_A.a$ for all $a\in A$, applying the bimodule maps in Definition~\ref{cfuadikf} gives $a.\alpha=\alpha.a$ and $a.\xi=\xi.a$. Next note that all the maps which are composed to give $\alpha$ are invertible, so the bimodule map is onto, and hence its image contains $1_A$. Therefore, there is an $a\in A$ with $a\,\alpha=1_A$. For the formula for $\alpha$, just write the evaluation in terms of the dual basis and use Proposition~\ref{xsrhkldjtyfggh}. Finally
\begin{eqnarray*}
\alpha^* &=& \<e^{i*},\overline{e^{j*}\,g_{ji}}\>  \cr
&=& \<e^{i*},g_{ji}^*\,\overline{e^{j*}}\>  \cr
&=& \<e^{i*},g_{ij}\,\overline{e^{j*}}\> \cr
&=& \<e^{i*}\,g_{ij},\overline{e^{j*}}\> \ .
\end{eqnarray*}
In terms of the dual basis, we have
\begin{eqnarray*}
\xi &=& (\star^{-1}G^{-1}\tens\id)(e_i\tens e^i) \cr
&=& (\star^{-1}\tens\id)(\overline{g_{ij}.e^j}\tens e^i)\cr
&=& e^{j*}\,g_{ji}\tens e^i\ .
\end{eqnarray*}
This formula, together with ${g_{ji}}^*=g_{ij}$, is enough to show that $\xi^*=\xi$. \quad$\square$

\medskip 
If the next corollary is worded somewhat strangely, it is because we were thinking about $C^*$-algebras and Hilbert $C^*$-modules, but wanted to make the statement rather more generally. 
This should be regarded as a comment, rather than an assumption that will prove vital.

\begin{cor} \label{vcjasfjjdghs}
Suppose that the matrix $g_{ij}\in M_n(A)$ can be factored into
a product of matrices with entries in $A$, $g_{ij}=r_{ik}\,r_{kj}$, so that the matrix $r_{ij}\in M_n(A)$ is also Hermitian. Then the elements $\alpha\in A$ and $\xi\in T_{\mathbb{N}}(L)$ are positive, in the sense that, for some $x_1,\dots,x_n\in L$,

\noindent
a)\quad $\alpha=\<x_1,\overline{x_1}\>+\dots+\<x_n,\overline{x_n}\>$ \ ,

\noindent
b)\quad $\xi=x_1\tens x_1^*+\dots+x_n\tens x_n^*$.
\end{cor}
\noindent {\bf Proof:}\quad We use the formulae given in Proposition~\ref{cfgyajkhj}. First
\begin{eqnarray*}
\alpha &=& \<e^{j*}\,g_{ji},\overline{e^{i*}}\>   \cr
&=&  \<e^{j*}\,r_{jk}\,r_{ki},\overline{e^{i*}}\>   \cr
&=&  \<e^{j*}\,r_{jk} , r_{ki}\,\overline{e^{i*}}\>   \cr
&=&  \<e^{j*}\,r_{jk} , \overline{e^{i*}\,r_{ki}^*}\>   \cr
&=&  \<e^{j*}\,r_{jk} , \overline{e^{i*}\,r_{ik}}\>   \ ,
\end{eqnarray*}
and then we put $x_k=e^{j*}\,r_{jk}$. Finally 
\begin{eqnarray*}
\xi &=& e^{j*}\,g_{ji}\tens e^i\cr
&=& e^{j*}\,r_{jk}\,r_{ki}\tens e^i\cr
&=& e^{j*}\,r_{jk}\tens r_{ki}\,e^i\cr
&=& e^{j*}\,r_{jk}\tens (e^{i*}\,r_{ki}^*)^*\cr
&=& e^{j*}\,r_{jk}\tens (e^{i*}\,r_{ik})^*\ .  \quad\square
\end{eqnarray*}

\begin{propos}\label{xgchjcjddf}
For all $e\in L$,
\begin{eqnarray*}
\xi\tens e\ =\ e\tens\xi.\Phi_L(\alpha)\in L\tens_A L\tens_A L\ ,
\end{eqnarray*}
where $\Phi_L: Z(A)\to Z(A)$ is the algebra map constructed in Proposition~\ref{cyjhgcxhj}.
\end{propos}
\noindent {\bf Proof:}\quad As $a.\xi=\xi.a$ for all $a\in A$, the map $e\mapsto \xi\tens e$ is a left module map from $L$ to $L\tens_A L\tens_A L$. It follows from Proposition~\ref{xdcxyjtalj} that there is an $\eta\in L\tens_A L$ so that $\xi\tens e=e\tens\eta\in L\tens_A L\tens_A L$ for all $e\in L$. Hence
\begin{eqnarray*}
(\ev(\id\tens G\,\star)\tens \id)(e\tens\eta) &=& 
\ev(\id\tens G\,\star)(\xi\tens e) \cr &=& \alpha.e\ =\ e.\Phi_L(\alpha)\ .
\end{eqnarray*}
As coevaluation is an isomorphism,  $(G\,\star\tens \id)\eta=\coev(\Phi_L(\alpha))$. The result follows by the definition of
$\xi$.\quad$\square$

\begin{cor} \label{xcvghjsfxsz}
We have $\alpha\,\Phi_L(\alpha)=1_A$. 
\end{cor}
\noindent {\bf Proof:}\quad Applying Proposition~\ref{xgchjcjddf} twice, for all $e\tens f\in L\tens_A L$,
\begin{eqnarray*}
\xi\tens e\tens f &=& e\tens\xi.\Phi_L(\alpha)\tens f  \in L\tens_A L\tens_A L\tens_A L \cr
&=&  e\tens\xi\tens f.\Phi_L(\Phi_L(\alpha))  \cr 
&=&  e\tens f\tens \xi.\Phi_L(\alpha)\,\Phi_L(\Phi_L(\alpha)) \ .
\end{eqnarray*}
But we can set (summing implicitly) $e\tens f=\xi$, giving
\begin{eqnarray*}
\xi\tens\xi\ =\ \xi\tens\xi.\Phi_L(\alpha)\,\Phi_L(\Phi_L(\alpha))\in L\tens_AL\tens_AL\tens_AL\ .
\end{eqnarray*}
This gives $\Phi_L(\alpha\,\Phi_L(\alpha))=1_A$ as required.\quad$\square$

\medskip As we shall see, it would be rather useful if $\alpha=1_A$. If $\Phi_L$ is the identity (as it is in the commutative case), then Corollary~\ref{xcvghjsfxsz} says that $\alpha^2=1_A$. As $\alpha$ is central and Hermitian, this is quite a restriction. If, in a subalgebra of a $C^*$-algebra,
$\alpha$ were positive (see Corollary~\ref{vcjasfjjdghs}), we would likely be able to recover $\alpha=1_A$ from this. But in general $\Phi_L$ is not the identity, and as it only depends on the bimodule isomorphism class of $L$, there is not much we can do about this. However there is one bit of Definition~\ref{cfuadikf} that we are reasonably free to change -- we can rescale the nondegenerate metric $G$.

\begin{propos}  \label{cvghajcds}
Two non-degenerate Hermitian metrics $G,G':\overline{L}\to L^\circ$ are related by $G'=R_z\, G$, where $R_z$ is right multiplication by an invertible central element $z\in Z(A)$. The corresponding inner products and the $\alpha$ (see Definition~\ref{cfuadikf}) are related by
\begin{eqnarray*}
\<e,\overline{f}\>' &=& \<e,\overline{f}\>\, z\ ,  \cr
\alpha' &=& \alpha\, \Phi_{L^\circ}(z^{-1}) \,z\ .
\end{eqnarray*}
The symmetry condition on the metric requires that $z$ is Hermitian, i.e.\ $z^*=z$.
\end{propos}
\noindent {\bf Proof:}\quad 
Recall that $G,G':\overline{L}\to L^\circ$ are invertible bimodule maps. Hence $G'\,G^{-1}:L^\circ\to L^\circ$ is an invertible bimodule map, and by using similar arguments to those used previously, it is right multiplication $R_z$ by an invertible central element $z\in Z(A)$. Now $G'=R_z\, G$, and substituting this into the expression in Definition~\ref{cfuadikf} gives
 $\alpha'\in A$ to be the image of $1_A$ under
\begin{eqnarray*}
A \stackrel{\coev} \longrightarrow L^\circ\tens_A L 
\stackrel{(G'\,\star)^{-1}\tens G'\,\star} \longrightarrow L\tens_A L^\circ 
\stackrel{\ev} \longrightarrow A\ .
\end{eqnarray*}
Since $\star$ is a bimodule map, this is the same as
\begin{eqnarray*}
A \stackrel{\coev} \longrightarrow L^\circ\tens_A L 
\stackrel{(G\,\star)^{-1}\tens G\,\star} \longrightarrow L\tens_A L^\circ 
\stackrel{ R_{z^{-1}} \tens R_z} \longrightarrow L\tens_A L^\circ 
\stackrel{\ev} \longrightarrow A\ .
\end{eqnarray*}
For an element $e\tens f\in L\tens_A L^\circ$ we have
\begin{eqnarray*}
\ev(R_{z^{-1}} \tens R_z)(e\tens f) &=& \ev(e\,z^{-1}\tens f\,z)  \cr
&=&  \ev(e\tens z^{-1}\, f) \,z \cr
&=& \ev(e\tens f)\,\Phi_{L^\circ}(z^{-1}) \,z\ .  \quad\square
\end{eqnarray*}

The conditions applied to the next proposition are again motivated by $C^*$-algebras, but restricting to $C^*$-algebras would be too strong, as it is quite likely that the result might be applied to smooth subalgebras of $C^*$-algebras or similar cases. 

\begin{propos}  \label{xfghjkjydxsx}
Suppose that $\alpha$ in Definition~\ref{cfuadikf} has a central Hermitian fourth root, i.e.\ 
$\beta=\beta^*\in Z(A)$ and $\beta^4=\alpha$. Rescale the metric $G$ (see Proposition~\ref{cvghajcds}) to get $G'=R_{\beta^{-2}}\,G$. Then $G'$ has corresponding $(\alpha')^2=1_A$. Also, if the original metric was positive, so is the new one, in the sense that
\begin{eqnarray*}
\<e,\overline{f}\>' &=& \beta^{-1*}\,\<e,\overline{f}\>\,\beta^{-1}\ .
\end{eqnarray*}
\end{propos}
\noindent {\bf Proof:}\quad  From Proposition~\ref{cvghajcds},
\begin{eqnarray*}
\alpha' &=& \alpha\, \Phi_{L^\circ}(\beta^2) \,\beta^{-2} \cr
&=& \beta^2\,\Phi_L^{-1}(\beta^2)\ =\ \Phi_L^{-1}(\beta^2\,\Phi_L(\beta^2))\ .
\end{eqnarray*}
Now Proposition~\ref{xcvghjsfxsz} shows that $(\alpha')^2=1_A$. \quad$\square$

\medskip In a $C^*$-algebra, if a positive element squares to the identity, then it must be the identity. 
This leads us to make the following definition, in the expectation that (from Corollary~\ref{vcjasfjjdghs} and Proposition~\ref{xfghjkjydxsx}) it should not be an uncommon possibility:

\begin{defin} \label{zdfgmhg}
A non-degenerate Hermitian metric $G$ is called {\em star compatible} if in Definition~\ref{cfuadikf} $\alpha=1$.
\end{defin}

\medskip
And here is where this account would likely end, were it not for  what is perhaps the most striking property of the element $\xi\in L^{\tens_A^2}$. 

\begin{propos} \label{xtdfsear} If the metric $G$ is star compatible, then 
any element of $L^{\tens_A^{2n}}$ can be written in the form $a.\xi^n$ for some $a\in A$. In particular, for $e,f\in L$,
\begin{eqnarray*}
e\tens f 
&=& \ev(e\tens (G\star)(f))\, .\, \xi\ .
\end{eqnarray*}
\end{propos}
\noindent {\bf Proof:}\quad Since $\xi$ commutes with all elements of $A$, it is enough to prove this for $n=1$, as repeated application of this will work for higher powers. For any $e\tens f\in L^{\tens_A^2}$, 
\begin{eqnarray}\label{vckfvcdytdss}
e\tens f &=& (\ev\tens\id^{\tens 2})(\id\tens\ev\tens\id^{\tens 3})(e\tens f\tens(\id\tens\coev\tens\id)\coev)\ .
\end{eqnarray}
As $\xi=((G\star)^{-1}\tens\id)\coev(1)$ commutes with elements of $L$ by Proposition~\ref{xgchjcjddf}, we can compute
\begin{eqnarray*}
(\id\tens\coev\tens\id)\coev(1)
&=& (\id\tens((G\star)(G\star)^{-1}\tens\id)\coev\tens\id)\coev(1)\cr
&=& (\id\tens(G\star)\tens\id^{\tens 2})(\id\tens\xi\tens\id)\coev(1)\cr
&=& (\id\tens(G\star)\tens\id^{\tens 2})(\coev(1)\tens\xi)\cr
&=& (\id\tens(G\star))\coev(1)\tens\xi\ .
\end{eqnarray*}
Substituting this into (\ref{vckfvcdytdss}) gives
\begin{eqnarray*}
e\tens f &=& \ev(\id\tens\ev\tens\id)\big(e\tens f\tens(\id\tens(G\star))\coev(1)\big)\,.\,\xi\cr
&=& \ev(e\tens (G\star)(f))\, .\, \xi\ .\quad\square
\end{eqnarray*}

\begin{theorem} \label{aehkjhcuy}
Suppose there is a star algebra $A$ and a left line module $L$ over $A$ with a star operation $\star:L\to \overline{L}$. In addition assume that there is a nondegenerate Hermitian metric $G: \overline{L}\to L^\circ$ which is star compatible. Then:

\noindent a)\quad The star algebra $T_{\mathbb{N}}(L)$ is isomorphic as an $\mathbb{N}$-graded algebra to
\begin{eqnarray*}
A\tens_{\mathbb{C}} P(\xi) \bigoplus L\tens_{\mathbb{C}} P(\xi) \ ,
\end{eqnarray*}
where $P(\xi)$ is the polynomial algebra in the variable $\xi$ of degree 2, and the degree of $L$ is $1$. The product is given by $\xi$ being central, the bimodule actions of $A$ on $L$, and $e\,f=\<e,\overline{f^*}\>\,\xi$ for $e,f\in L$.The star operation is given by that specified on $A$ and $L$, together with $\xi^*=\xi$. 

\noindent b)\quad The even-odd graded algebra $A\oplus L$ (with product $e\,f=\<e,\overline{f^*}\>$ for $e,f\in L$) is a $\mathbb{Z}/2$-graded Hopf-Galois extension of $A$. 
\end{theorem}
\noindent {\bf Proof:}\quad Combining the results in this section gives (a). To check consistency of the odd-odd product with the star operation, we need to check that
\begin{eqnarray*}
\ev(e\tens (G\star)(f))^* \ =\ \ev(f^*\tens (G\star)(e^*))\ ,
\end{eqnarray*}
but this simply becomes $\<e,\overline{f^*}\>^*=\<f^*,\overline{e}\>$, the usual symmetry relation for the inner product. To get (b), just quotient by the relation $\xi=1$. \quad$\square$

\medskip We should comment that Theorem~\ref{aehkjhcuy} is just what we should expect from the classical theory. Consider a locally trivial real line bundle over a compact Hausdorff topological space $X$. The star algebra $T_{\mathbb{N}}(L)$ corresponds to the functions on the line bundle
which are polynomial in the fibre direction. The grading  of $T_{\mathbb{N}}(L)$ corresponds to the order of the polynomial.  

The transition functions for the line bundle can be taken to be in $\mathbb{R}^*$. However if there is a metric on the bundle, we can reduce the transition functions to have values $\pm 1$. If we consider functions which are even in the fibre direction (i.e.\ $f(-x)=f(x)$), then the $\pm 1$ makes no difference, and we may as well look at even functions on the trivial bundle $X\times \mathbb{R}$. Now the odd functions are given by some (not identically zero) linear odd function times even functions. (Of course, this requires a regularity condition, but we are not at present considering all continuous functions.) However if the bundle is not trivial, we cannot continuously choose a global non-vanishing linear function on all the fibres, so we are stuck with several choices on several open sets. 

Again using the metric, we can restrict to the points on the line bundle which are distance one from the zero section. This is a double cover of $X$, and functions on this space correspond to the 
$\mathbb{Z}/2$-graded Hopf-Galois extension of $A$ mentioned in Theorem~\ref{aehkjhcuy}.

\section{Representations of functions on $\mathbb{R}$} \label{cvaufv1}
Here we will generalise the algebra in Theorem~\ref{aehkjhcuy} to functions other than polynomials.
In particular we assume the conditions for Theorem~\ref{aehkjhcuy}: 
That is we assume that there is a star algebra $A$ and a left line module $L$ over $A$ with a star operation $\star:L\to \overline{L}$. In addition assume that there is a nondegenerate Hermitian metric $G: \overline{L}\to L^\circ$ which is star compatible.
We use $x$ as the standard coordinate function on $\mathbb{R}$.
We identify $\xi\in L\tens_A L$ with the function $x^2$ on $\mathbb{R}$. Generically we have
a $\mathbb{Z}/2$-graded algebra
\begin{eqnarray} \label{bchsdlvcyy}
B_L \ =\ 
\big\{(f_0,f_1) \,|\, f_0:\mathbb{R}\to A\ ,\ \ f_1:\mathbb{R}\to L\big\}\ .
\end{eqnarray}
We have not yet specified exactly which classes of functions are to be used in 
(\ref{bchsdlvcyy}), but the product is given by $(f_0,f_1)(g_0,g_1)=(h_0,h_1)$, where
\begin{eqnarray} \label{qbchsdlvcty}
h_0(x) &=& f_0(x)\, g_0(x) + x^2\,\<f_1(x),\overline{g_1(x)^*}\>\ ,\cr
h_1(x) &=& f_0(x)\, g_1(x) +f_1(x)\, g_0(x) \ .
\end{eqnarray}
This product is not at all random -- it is chosen to generalise the case in Section~\ref{cvaufv2}
where the functions are polynomials in the fibre direction. However we still need to check associativity:

\begin{propos} \label{vcagjkscdrs}
 In the star compatible case, the product in (\ref{qbchsdlvcty}) and the star operation 
\[(f_0(x),f_1(x))^*=(f_0(x)^*,f_1(x)^*)\]
 make $B_L$ into an associative star algebra.
\end{propos}
\noindent {\bf Proof:}\quad We need to show that
\begin{eqnarray*}
&&\big((f_0(x),f_1(x))\,(g_0(x),g_1(x))\big)\,(h_0(x),h_1(x))\cr &=& (f_0(x),f_1(x))\,\big((g_0(x),g_1(x))\,(h_0(x),h_1(x))\big) \ .
\end{eqnarray*}
The products with all 0 indices are just products in $A$, so the associative law holds. The products with one 1 index are associative just by the usual properties of bimodules. The products with two 1 indices are:
\begin{eqnarray*}
f_0(x)\,\big(g_1(x)\,h_1(x)\big) &=& x^2\,f_0(x)\,\<g_1(x),\overline{h_1(x)^*}\> \cr
&=& x^2\,\<f_0(x)\,g_1(x),\overline{h_1(x)^*}\>\cr
&=& \big(f_0(x)\,g_1(x)\big)\,h_1(x)\ ,\cr
f_1(x)\,\big(g_0(x)\,h_1(x)\big) &=& x^2\,\<f_1(x),\overline{h_1(x)^*\,g_0(x)^*}\> \cr
&=& x^2\,\<f_1(x),g_0(x)\,\overline{h_1(x)^*}\> \cr
&=& x^2\,\<f_1(x)\,g_0(x),\overline{h_1(x)^*}\> \cr
&=& \big(f_1(x)\,g_0(x)\big)\,h_1(x)\ ,\cr
f_1(x)\,\big(g_1(x)\,h_0(x)\big) &=&  x^2\,\<f_1(x),\overline{h_0(x)^*\,g_1(x)^*}\> \cr
&=&  x^2\,\<f_1(x),\overline{g_1(x)^*}\,h_0(x)\> \cr
&=&  x^2\,\<f_1(x),\overline{g_1(x)^*}\>\,h_0(x) \cr
&=& \big(f_1(x)\,g_1(x)\big)\,h_0(x)\ .
\end{eqnarray*}
The most difficult case is three 1 indices:
\begin{eqnarray*}
f_1(x)\,\big(g_1(x)\,h_1(x)\big) &=&  x^2\,f_1(x)\,\<g_1(x),\overline{h_1(x)^*}\> \ ,\cr
\big(f_1(x)\,g_1(x)\big)\,h_1(x) &=&  x^2\,\<f_1(x),\overline{g_1(x)^*}\>\,h_1(x)\ .
\end{eqnarray*}
To show that the right hand sides are mutually equal we need to verify that
\begin{eqnarray}  \label{vcyuafxds}
(\id\tens\ev)(\id^{\tens 2}\tens G\,\star) \ =\ \ev(\id\tens G\,\star)\tens\id:L\tens_A L\tens_A L\to L\ .
\end{eqnarray}
By Definitions~\ref{zdfgmhg} and \ref{cfuadikf}, 
\begin{eqnarray*}
\ev(\id\tens G\,\star) \ =\ \coev^{-1}(G\,\star\tens\id):L\tens_A L\to A\ ,
\end{eqnarray*}
so (\ref{vcyuafxds}) becomes 
\begin{eqnarray*}
\id\tens\coev^{-1}\ =\ \ev\tens\id:L\tens_A L^\circ\tens_A L\to L\ ,
\end{eqnarray*}
and these are shown to be equal in Proposition~\ref{bhcsavch}. 

To check the star algebra property, we need to verify
\begin{eqnarray*}
(g_0(x)^*,g_1(x)^*)\,(f_0(x)^*,f_1(x)^*) &=& \big((f_0(x),f_1(x))\,(g_0(x),g_1(x))\big)^*\ .
\end{eqnarray*}
This is trivial, apart from the two 1 index product. In this case 
\begin{eqnarray*}
\big(g_1(x)^*\, f_1(x)^*\big)^* &=& \<g_1(x)^*,\overline{f_1(x)}\>^*\cr
&=& \<f_1(x),\overline{g_1(x)^*}\> \cr
&=& f_1(x)\, g_1(x)\ ,
\end{eqnarray*}
as required. $\quad\square$

\medskip Next we specify a vector space for $B_L$ to act on. If we consider $A\subset \mathcal{B}(\mathcal{H})$ (the bounded linear operators on the Hilbert space $\mathcal{H}$), we define a new vector space as a sum $\mathcal{H}_0\oplus \mathcal{H}_1$, where
\begin{eqnarray}  \label{cfghjhgcz}
\mathcal{H}_0 \ =\  C(\mathbb{R}_\infty,\mathcal{H})\ ,\ \ 
\mathcal{H}_1 \ =\  C_c(\mathbb{R},\overline{L}\tens_A\mathcal{H})\ .
\end{eqnarray}
Here $C(\mathbb{R}_\infty,\mathcal{H})$ is simply the algebra of continuous functions from $\mathbb{R}_\infty$ (the one point compactification of $\mathbb{R}$, topologically a circle) to $\mathcal{H}$. Then $C_c(\mathbb{R},\overline{L}\tens_A\mathcal{H})$ is a sum of $\overline{e^i}\tens k_i(x)$, where the $e^i$ are elements of the basis of $L$ (see Definition~\ref{canonyy}) and the $k_i:\mathbb{R}\to \mathcal{H}$ are continuous functions of compact support.

We feel obliged to give an apology for the way the bars are going to work out -- this is due to the convention of writing a $C^*$-algebra left acting on a Hilbert space.
If we write the inner product on $\mathcal{H}$ as $\<,\>_\mathcal{H}:\overline{\mathcal{H}} \tens \mathcal{H}  \to \mathbb{C}$ and the Hilbert $C^*$-module inner product on $L$ as $\<,\>$, then we define (generalised) inner products $\<,\>_0$ on 
$\mathcal{H}_0$ and $\<,\>_1$ on $\mathcal{H}_1$ as
the following continuous $\mathbb{C}$ valued functions on $\mathbb{R}$, respectively
\begin{eqnarray} \label{mkiobvycyudh}
\<\overline{v_0},w_0\>_0 &=& x\mapsto \<\overline{v_0(x)},w_0(x)\>_\mathcal{H}\ ,\cr
\<\overline{v_1},w_1\>_1 &=& x\mapsto  
\<,\>_\mathcal{H}(\id\tens \<,\>\la\id) ((\id\tens\mathrm{bb}^{-1})\Upsilon(\overline{v_1(x)})\tens w_1(x)) \ .
\end{eqnarray}  
The first equation in (\ref{mkiobvycyudh}) is quite simple -- evaluate the functions on $\mathbb{R}$ pointwise and apply the Hilbert space inner product. The second is rather more complicated. To explain it, write $v_1(x)=\overline{e}\tens u$ and $w_1(x)=\overline{e'}\tens u'$, and then, where $\la$ is used to denote the action of $A$ on $\mathcal{H}$,
\begin{eqnarray}  \label{vcuaykhjgxzaz}
&& 
\<,\>_\mathcal{H}(\id\tens \<,\>\la\id) ((\id\tens\mathrm{bb}^{-1})\Upsilon(\overline{v_L(x)})\tens w_L(x))  \cr
&=&  \<,\>_\mathcal{H}(\id\tens \<,\>\la\id) ((\id\tens\mathrm{bb}^{-1})\Upsilon(\overline{\overline{e}\tens u})\tens \overline{e'}\tens u') \cr
&=&  \<,\>_\mathcal{H}(\id\tens \<,\>\la\id) ((\id\tens\mathrm{bb}^{-1})(\overline{u} \tens \overline{\overline{e}} )\tens \overline{e'}\tens u') \cr  
&=&  \<,\>_\mathcal{H}(\id\tens \<,\>\la\id) (\overline{u} \tens e \tens \overline{e'}\tens u') \cr 
&=&  \<,\>_\mathcal{H} (\overline{u} \tens \<e , \overline{e'}\>\la u') \cr 
&=&  \< \overline{u} , \<e , \overline{e'}\>\la u'\>_\mathcal{H}\ .
\end{eqnarray}

\begin{propos}  \label{xzerytttt}
The formula $(f_0,f_1)\la (v_0,v_1)=(w_0,w_1)$, where
\begin{eqnarray*}
w_0(x) &=& f_0(x)\la\, v_0(x)+x^2\,(\<,\>\la\,\id)(f_1(x)\tens v_1(x))\ ,\cr
w_1(x) &=& \overline{f_1(x)^*}\tens_A v_0(x) + f_0(x)\la\, v_1(x)\ ,
\end{eqnarray*}
gives an action of the algebra $B_L$ (see (\ref{bchsdlvcyy})) on $\mathcal{H}_0\oplus \mathcal{H}_1$. Note that if we write $v_1(x)=\overline{e}\tens u$, then the two less obvious terms above are
\begin{eqnarray*}
(\<,\>\la\,\id)(f_1(x)\tens v_1(x)) &=& \< f_1(x) , \overline{e}\>\la\, u\ ,  \cr
f_0(x)\la\, v_1(x) &=& f_0(x)\la\, \overline{e}\tens u \cr
&=&  \overline{e\,\ra\, f_0(x)^*}\tens u  \ .
\end{eqnarray*}
\end{propos}
\noindent {\bf Proof:}\quad 
Set $(g_0,g_1)\la (w_0,w_1)=(y_0,y_1)$, where
\begin{eqnarray*}
y_0(x) &=& (g_0(x)\,f_0(x))\la\, v_0(x)+x^2\,g_0(x)\< f_1(x) , \overline{e}\>\la\, u \cr
&& +\ x^2\,\<g_1(x),\overline{f_1(x)^*}\>\,\la\, v_0(x) 
+ x^2\,\<g_1(x),f_0(x)\la\, \overline{e}\>\,\la\, u \cr
&=& \big(g_0(x)\,f_0(x) + x^2\,\<g_1(x),\overline{f_1(x)^*}\>\big)\la\, v_0(x) \cr
&&+\ x^2\,\<\big(g_1(x)\,f_0(x) + g_0(x)\,f_1(x)  \big),\overline{e}\>\,\la\, u \ ,
\end{eqnarray*}
as required. Similarly, where $w_1(x)=\overline{e'}\tens u'$,
\begin{eqnarray*}  \label{vchjakvc}
y_1(x) &=& \overline{g_1(x)^*}\tens_A w_0(x) +  g_0(x) \, \la\, \overline{e'}\tens u' \cr
&=& \overline{g_1(x)^*}\tens_A f_0(x)\la\, v_0(x) + x^2\,\overline{g_1(x)^*}\tens_A 
\< f_1(x) , \overline{e}\>\la\, u \cr
&&  +\  g_0(x)\la   \overline{f_1(x)^*}\tens_A v_0(x) + (g_0(x)\,f_0(x))\la\, v_1(x) \cr
&=& \overline{(g_1(x)\,f_0(x))^*}\tens_A v_0(x) + x^2\,\overline{g_1(x)^*}\tens_A 
\< f_1(x) , \overline{e}\>\la\, u \cr
 && +\    \overline{(g_0(x)\,f_1(x))^*}\tens_A v_0(x) + (g_0(x)\,f_0(x))\la\, v_1(x) 
\end{eqnarray*}
and this is not so obvious. We would require the following equality
to get an action:
\begin{eqnarray*}
\<g_1(x),\overline{f_1(x)^*}\> \, \la\,\overline{e}\tens u &=& \overline{g_1(x)^*}\tens_A 
\< f_1(x) , \overline{e}\>\la\, u\ ,
\end{eqnarray*}
which can be simplified to showing
\begin{eqnarray*}
\<g_1(x),\overline{f_1(x)^*}\> \, \overline{e} &=& \overline{g_1(x)^*}\,
\< f_1(x) , \overline{e}\>\ .
\end{eqnarray*}
This is implied by the following equality, proved in Lemma~\ref{xzeryt},
\begin{eqnarray*}
(\ev\tens\id)(\id\tens G\,\star\tens\id) &=& (\id\tens\ev)(\star\tens \id\tens G)
:L\tens_A L\tens_A\overline{L}\to \overline{L}\ .\quad\square
\end{eqnarray*}

\begin{lemma}  \label{xzeryt}
Assuming that the metric is star compatible, 
\begin{eqnarray*}
(\ev\tens\id)(\id\tens G\,\star\tens\id) &=& (\id\tens\ev)(\star\tens \id\tens G)
:L\tens_A L\tens_A\overline{L}\to \overline{L}\ .
\end{eqnarray*}
\end{lemma}
\noindent {\bf Proof:}\quad This is equivalent to showing
\begin{eqnarray*}
(\ev\tens\id)(\id\tens G\,\star\tens\star) &=& (\id\tens\ev)(\star\tens \id\tens G\,\star)
:L\tens_A L\tens_A L\to \overline{L}\ . 
\end{eqnarray*}
As these are both bimodule maps, it is enough to prove that these maps coincide on $e\tens\xi\in L^{\tens_A^3}$. As $\xi$ is central, then $e\tens\xi=\xi\tens e$, so it is enough to prove that 
\begin{eqnarray}   \label{xfghjcfhj}
(\ev\tens\id)(\id\tens G\,\star\tens\star)(\xi\tens e) &=& (\id\tens\ev)(\star\tens \id\tens G\,\star)(e\tens\xi)\ .
\end{eqnarray}
Note that $\ev(\id\tens G\,\star)\xi=1_A$ by star compatibility, thus
verifying (\ref{xfghjcfhj}).\quad   $\square$

\section{$C^*$-algebra completions and the Thom construction}  \label{cvaufv3}
We suppose that $A$ is a unital $C^*$-algebra
with norm written $|a|_A$ for $a\in A$, and that $L$ is a left line module for $A$. 
We take a nondegenerate positive inner product 
 $\<,\>:L\tens_A \overline{L}\to A$ which is a Hilbert $C^*$-module on
$L$. The corresponding matrix $g^{ij}\in M_n(A)$ (see Proposition~\ref{xsrhkldjtyfggh}) we take to have Hermitian square root $r^{ij}$, i.e.\ 
$r^{ip}\,r^{pj}=g^{ij}$. We also assume that there is a constant $M\ge 0$ so that, for all $i,j$
and all $a\in A$,
\begin{eqnarray}  \label{vcxghajsxr}
\big|e_i(e^j.a)|_A\ \le\ M\,|a|_A\ .
\end{eqnarray}
The direct sum $\mathcal{H}_0\oplus \mathcal{H}_1$ is made into a Hilbert space as follows:

\begin{defin}
Define  a $\mathbb{C}$ valued inner product on $\mathcal{H}_0\oplus \mathcal{H}_1$ 
in terms of the function valued ones from (\ref{mkiobvycyudh}) by 
\begin{eqnarray} \label{poivgcgh}
\<\overline{(v_0,v_1)},(w_0,w_1)\>_{0,1} &=& \int_{\mathbb{R}} 
\Big( \frac{\<\overline{v_0},w_0\>_0 +  x^2\,\<\overline{v_1},w_1\>_1}{1+x^2}
\Big)\, \extd x\ .
\end{eqnarray}
The definition of $\mathcal{H}_0$ and $\mathcal{H}_1$ in (\ref{cfghjhgcz}) guarantees that this is
finite, and positivity is guaranteed by the definitions of the inner products in 
(\ref{mkiobvycyudh}). Now define $\mathcal{H}_{0,1}$ to be the Hilbert space given by completing
$\mathcal{H}_0\oplus \mathcal{H}_1$ with inner product $\<,\>_{0,1}$. 
\end{defin}

 If the readers have any doubts over the positivity of $\<,\>_{1}$, Lemma~\ref{kjydfdghj}
should convince them. Note that we spell out the summation over the indices explicitly in Lemma~\ref{kjydfdghj}. In future, we will assume that repeated indices inside a norm are summed before the norm is taken, as in this case. 

\begin{lemma}  \label{kjydfdghj}
For any $v_1(x)=\sum_i\overline{e^i}\tens k_i(x)\in \mathcal{H}_1$, 
\begin{eqnarray*}
\<\overline{v_1(x)},v_1(x)\>_1 &=& \sum_p \Big\| \sum_j r^{pj}\,k_j(x) \Big\|_\mathcal{H}^2 \ .
\end{eqnarray*}
\end{lemma}
\noindent {\bf Proof:}\quad From (\ref{vcuaykhjgxzaz}),
\begin{eqnarray*}
\<\overline{v_1(x)},v_1(x)\>_1 &=& \<\overline{k_i(x)}\,,\, \<e^i,\overline{e^j}\>\, k_j(x)\>_\mathcal{H} \cr
&=&  \<\overline{k_i(x)}\,,\, g^{ij}\, k_j(x)\>_\mathcal{H} \cr
&=&  \<\overline{r^{pi}\,k_i(x)}\,,\, r^{pj}\, k_j(x)\>_\mathcal{H} \ .\quad\square
\end{eqnarray*}

\begin{lemma}  \label{kjydfdghtyj1}
There is a constant $C_1\ge 0$ (depending only on $A$, $L$, the dual basis of $L$
and the Hilbert $C^*$-module structure on $L$) so that
\begin{eqnarray*}
\<\overline{f_0(x)\la\, v_1(x)},f_0(x)\la\, v_1(x)\>_1 &\le& C_1\,|f_0(x)|_A^2\,\<\overline{v_1(x)},v_1(x)\>_1\ .
\end{eqnarray*}
\end{lemma}
\noindent {\bf Proof:}\quad If we write $v_1(x)=\overline{e^i}\tens k_i(x)$ (summation implicit), then
\begin{eqnarray*}
f_0(x)\la\, v_1(x) &=& f_0(x)\la\, \overline{e^i}\tens k_i(x) \cr
&=&  \overline{e^i\,\ra\, f_0(x)^*}\tens k_i(x)  \cr
&=&  \overline{e_j(e^i\, f_0(x)^*)\,e^j}\tens k_i(x)   \cr
&=&  \overline{e^j}\tens e_j(e^i\, f_0(x)^*)^*\,k_i(x) \ .
\end{eqnarray*}
By Lemma~\ref{kjydfdghj}, 
\begin{eqnarray*}
\<\overline{f_0(x)\la\, v_1(x)},f_0(x)\la\, v_1(x)\>_1
&=&  \sum_p \big\|  r^{pj}\,e_j(e^i\, f_0(x)^*)^*\,k_i(x) \big\|_\mathcal{H}^2 \ .
\end{eqnarray*}
Using the dual basis property and Proposition~\ref{xsrhkldjtyfggh} we can write
\begin{eqnarray*}
r^{pj}\,e_j(e^i\, f_0(x)^*)^*\,k_i(x) &=& r^{pj}\,e_j(e^m\, f_0(x)^*)^*\,e_m(e^i)^*\,k_i(x) \cr
&=&  r^{pj}\,e_j(e^m\, f_0(x)^*)^*\,g_{mq}\,g^{qi}\,k_i(x) \cr
&=&  r^{pj}\,e_j(e^m\, f_0(x)^*)^*\,g_{mq}\,r^{qs}\,r^{si}\,k_i(x) \ ,\cr
\big\|r^{pj}\,e_j(e^i\, f_0(x)^*)^*\,k_i(x) \big\|_\mathcal{H} &\le & 
\big|r^{pj}\,e_j(e^m\, f_0(x)^*)^*\,g_{mq}\,r^{qs}\big|_A     \ 
\big\|  r^{si}\,k_i(x)     \big\|_\mathcal{H} \ .
\end{eqnarray*}
We use the Cauchy-Schwarz inequality on this sum over $s$ to obtain
\begin{eqnarray*}
&& \big\|r^{pj}\,e_j(e^i\, f_0(x)^*)^*\,k_i(x) \big\|_\mathcal{H}^2  \cr
&\le & \Big(\sum_s  \big|r^{pj}\,e_j(e^m\, f_0(x)^*)^*\,g_{mq}\,r^{qs}\big|_A^2   \Big)\  
\Big(\sum_s    \big\|  r^{si}\,k_i(x)     \big\|_\mathcal{H}^2 \Big)  \cr
&\le & \Big(\sum_s  \big|r^{pj}\,e_j(e^m\, f_0(x)^*)^*\,g_{mq}\,r^{qs}\big|_A^2   \Big)\  
\<\overline{v_1(x)},v_1(x)\>_1\ , 
\end{eqnarray*}
and summing this over $p$ gives
\[
\<\overline{f_0(x)\la\, v_1(x)},f_0(x)\la\, v_1(x)\>_1 
\le  \Big(\sum_{s,p}  \big|r^{pj}\,e_j(e^m\, f_0(x)^*)^*\,g_{mq}\,r^{qs}\big|_A^2   \Big)\  
\<\overline{v_1(x)},v_1(x)\>_1\ . 
\]
The result follows from (\ref{vcxghajsxr}).\quad$\square$

\begin{lemma}   \label{kjydfdghtyj2}
\[
\big\<\overline{(\<,\>\la\,\id)(f_1(x)\tens v_1(x))},(\<,\>\la\,\id)(f_1(x)\tens v_1(x))\big\>_0
\le \big| \< f_1(x) , \overline{f_1(x)}\>\big|_A\  \<  \overline{ v_1(x) }, v_1(x)\>_1  \ .
\]
\end{lemma}
\noindent {\bf Proof:}\quad Write $v_1(x)=\overline{e^i}\tens k_i(x)$ (summation implicit), so that
\begin{eqnarray*}
(\<,\>\la\,\id)(f_1(x)\tens v_1(x)) &=& \< f_1(x) , \overline{e^i}\>\la\, k_i(x)\ .
\end{eqnarray*}
Then
\begin{eqnarray}  \label{aethjtyd}
&&\big\<\overline{(\<,\>\la\,\id)(f_1(x)\tens v_1(x))},(\<,\>\la\,\id)(f_1(x)\tens v_1(x))\big\>_0\cr
&=& \big\<  \overline{ \< f_1(x) , \overline{e^i}\>\, k_i(x) }, \<f_1(x) , \overline{e^i}\>\, k_i(x) \big\>_\mathcal{H}  \cr
&=& \big\<  \overline{ k_i(x) }, \< e^i, \overline{ f_1(x)}\>\, \<f_1(x) , \overline{e^i}\>\, k_i(x) \big\>_\mathcal{H}  \ .
\end{eqnarray}
Recall from \cite{Lance}
(with a brief check that the different side used for the conjugate does not matter) that for a Hilbert $C^*$-module there is a version of the Cauchy-Schwartz lemma as follows, in terms of inequalities of positive operators
\begin{eqnarray*}
\<y,\overline{x}\>\, \<x,\overline{y}\> \ \le\ |\<x,\overline{x}\>|_A\ \<y,\overline{y}\>\ ,
\end{eqnarray*}
and using this in (\ref{aethjtyd}) gives
\begin{eqnarray*} 
&&
\big\<\overline{(\<,\>\la\,\id)(f_1(x)\tens v_1(x))},(\<,\>\la\,\id)(f_1(x)\tens v_1(x))\big\>_0\cr
&\le& 
\big| \< f_1(x) , \overline{f_1(x)}\>\big|_A\  \big\<  \overline{ k_i(x) }, \<e^i , \overline{e^i}\>\, k_i(x) \big\>_\mathcal{H}  \ .
\end{eqnarray*}
Finally, use the definition of $\<,\>_1$ again.\quad$\square$

\begin{propos}  \label{vcahjgkvyfdf1}
Elements $(f_0,f_1)\in B_L$ act on $\mathcal{H}_{0,1}$  (with
inner product given by (\ref{poivgcgh})) as linear operators, with operator norm bounded by
a constant (depending only on $A$, $L$, the dual basis of $L$
and the Hilbert $C^*$-module structure on $L$) times the square root of
\begin{eqnarray*}
\sup_{x\in\mathbb{R}} |f_0(x)|_A^2+\sup_{x\in\mathbb{R}} x^2\,\Big(
 |\<f_1(x)^* , \overline{f_1(x)^*}\>|_A+ |\< f_1(x) , \overline{f_1(x)}\>\big|_A\Big)\ .
\end{eqnarray*}
\end{propos}
\noindent {\bf Proof:}\quad We check this on the dense subset $\mathcal{H}_{0}\oplus \mathcal{H}_{1}$.  Set $(f_0,f_1)\la (v_0,v_1)=(w_0,w_1)$, where
\begin{eqnarray} \label{cabahcvyy}
w_0(x) &=& f_0(x)\la\, v_0(x)+x^2\,(\<,\>\la\,\id)(f_1(x)\tens v_1(x))\ ,\cr
w_1(x) &=& \overline{f_1(x)^*}\tens_A v_0(x) + f_0(x)\la\, v_1(x)\ .
\end{eqnarray}
The two less obvious terms in 
(\ref{cabahcvyy}) (those involving $v_1(x)$) are dealt with in Lemmas \ref{kjydfdghtyj1} and \ref{kjydfdghtyj2}, and we will import these results. For the other two terms,
(using $\|.\|_0^2$ and $\|.\|_1^2$ for the $\<,\>_0$ and $\<,\>_1$ inner product of an element with itself),
\begin{eqnarray*}
\|f_0(x)\la\, v_0(x)\|_0  &\le &  |f_0(x)|_A\, \|v_0(x)\|_\mathcal{H}\ ,\cr
\| \overline{f_1(x)^*}\tens_A v_0(x) \|_1^2 &=&
 \< \overline{v_0(x) } , \<f_1(x)^* , \overline{f_1(x)^*}\>\la v_0(x) \>_\mathcal{H} \cr
 &\le & |\<f_1(x)^* , \overline{f_1(x)^*}\>|_A\ \|v_0(x) \|_\mathcal{H} ^2\ .
\end{eqnarray*}
Now we write out the contributions to $\<\overline{(w_0,w_1)},(w_0,w_1)\>_{0,1}$ as follows:
\begin{eqnarray*}
\int_{\mathbb{R}} \frac{\|f_0(x)\la\, v_0(x)\|_0^2}{1+x^2}\ \extd x &\le & 
\int_{\mathbb{R}} \frac{|f_0(x)|_A^2\, \|v_0(x)\|_0^2}{1+x^2}\ \extd x 
\cr
\int_{\mathbb{R}} \frac{\|x^2\,(\<,\>\la\,\id)(f_1(x)\tens v_1(x))\|_0^2}{1+x^2}\ \extd x &\le & 
\int_{\mathbb{R}} \frac{x^4\,\big| \< f_1(x) , \overline{f_1(x)}\>\big|_A^2\, \|v_1(x)\|_1^2}{1+x^2}\ \extd x 
\cr
\int_{\mathbb{R}}   \frac{x^2}{1+x^2} \,  \| \overline{f_1(x)^*}\tens_A v_0(x) \|_1^2\ \extd x &\le & 
\int_{\mathbb{R}}  \frac{x^2}{1+x^2} \,  |\<f_1(x)^* , \overline{f_1(x)^*}\>|_A\ \|v_0(x) \|_0^2\ \extd x 
\cr
\int_{\mathbb{R}}  \frac{x^2}{1+x^2}   \, \| f_0(x)\la\, v_1(x) \|_1^2\ \extd x &\le & C_1\,
\int_{\mathbb{R}}   \frac{x^2}{1+x^2}  \, |f_0(x)|_A^2\,\ \|v_1(x) \|_1^2\ \extd x \ .\quad\square
\end{eqnarray*}

\medskip The adjoint of a bounded operator $T:\mathcal{H}_{0,1}\to\mathcal{H}_{0,1}$ is an operator $T^*:\mathcal{H}_{0,1}\to\mathcal{H}_{0,1}$ defined so that
\begin{eqnarray*}
\<\overline{(u_0,u_1)},T(v_0,v_1)\>_{0,1} &=& \<\overline{T^*(u_0,u_1)},(v_0,v_1)\>_{0,1}\ .
\end{eqnarray*}

\begin{propos}   \label{vcahjgkvyfdf2}
The star operation on $B_L$ (see Proposition~\ref{vcagjkscdrs}) gives a star representation of the algebra on bounded operators on 
$\mathcal{H}_{0,1}$. 
\end{propos}
\noindent {\bf Proof:}\quad Set $v_1(x)=\overline{e^i}\tens k_i(x)$ and $u_1(x)=\overline{e^i}\tens s_i(x)$. There are four terms to check to verify that the star operations coincide, the easiest one being
\begin{eqnarray*}
\<\overline{u_0(x)},f_0(x).v_0(x)\>_0 &=& \<\overline{f_0(x)^*.u_0(x)},v_0(x)\>_0\ .
\end{eqnarray*}
Next consider
\begin{eqnarray*}
&& \<\overline{u_0(x)},(\<,\>\tens\id)(f_1(x)\tens v_1(x))\>_0 \cr
&=& \<\overline{u_0(x)},\<f_i(x),\overline{e^i}\>\,k_i(x)\>_\mathcal{H} \cr
&=& \big\< \overline{\overline{f_1(x)}\tens u_0(x)},
\overline{e^i}\tens k_i(x)\big\>_1\ .
\end{eqnarray*}
Then
\begin{eqnarray*}
\<\overline{u_1(x)},f_0(x)\,\la\,v_1(x))\>_1 &=& \big\< \overline{s_j(x)},\<e^j,f_0(x)\,
\overline{e^i}\>\,k_i(x)\big\>_\mathcal{H} \cr
&=& \big\< \overline{s_j(x)},\<e^j\,f_0(x),
\overline{e^i}\>\,k_i(x)\big\>_\mathcal{H} \cr
&=& \big\< \overline{f_0(x)^*\,\overline{e^j}\tens s_j(x)},\overline{e^i}\tens k_i(x)\big\>_1\cr
  &=& \big\< \overline{f_0(x)^*\,\la\,u_1(x)},v_1(x)\big\>_1\ .
\end{eqnarray*}
Finally
\begin{eqnarray*}
\<\overline{u_1(x)},\overline{f_1(x)^*}\tens v_0(x))\>_1 &=&
 \big\< \overline{s_j(x)},\<e^j,
\overline{f_1(x)^*}\>\,v_0(x)\big\>_\mathcal{H} \cr
&=&  \big\< \overline{\<e^j,
\overline{f_1(x)^*}\>^*\,s_j(x)},v_0(x)\big\>_\mathcal{H} \cr
&=&  \big\< \overline{\<f_1(x)^*,\overline{e^j}\>\,s_j(x)},v_0(x)\big\>_\mathcal{H} \ .\quad\square
\end{eqnarray*}

\medskip
Suppose that $X$ is a compact topological space. Given an $\mathbb{R}^n$ bundle on $X$, the Thom construction adds a single point at infinity common to all the fibres, so it is just the one point compactification of the total space of the bundle. Alternatively, we could add a point at infinity to each fibre separately (i.e.\ we add a whole copy of $X$), so that we get an asscoiated $S^n$ bundle. 
For noncommutative line bundles (the $n=1$ case) and unital $C^*$-algebras $A$, we can perform both of these compactifications, to get unital $C^*$-algebras. 

One question should be raised now, before it causes confusion.
A circle bundle from a line bundle?  Surely we did this by the $\mathbb{Z}$-graded algebra in Section \ref{vcuskftsd}, with star structure included in Section \ref{vuwicvhcv}? 
In Section \ref{cytjadsdrta} it was shown that there was a $\mathbb{CZ}$ coaction on 
the $\mathbb{Z}$-graded algebra, and this should give a circle group coaction. By now the reader should be becoming suspicious, there is no reason why a circle bundle associated to a line bundle should be a circle group principal bundle. The confusion is caused by taking the wrong star structure. The star structure given by functions on the total space of a classical real line bundle follows the pattern of Section \ref{xctkvdft}, and the relevant algebra is the $\mathbb{N}$-graded algebra.

\begin{theorem}
From the action of $B_L$ on the Hilbert space $\mathcal{H}_{0,1}$  three $C^*$-algebras can be formed. The functions $(f_0,f_1)\in B_L$ have $f_1:\mathbb{R}\to L$ chosen so that
\begin{eqnarray*}
f_1(x)\ =\ \sum_i f_{1,i}(x)\,e^i\ ,
\end{eqnarray*}
where each $f_{1,i}:\mathbb{R}\to A$ is of compact support. 
The choice of $f_0:\mathbb{R}\to A$ depends on the case, as given below. Then the $C^*$ completion is taken as operators on $\mathcal{H}_{0,1}$. 

\noindent a)\quad The non-unital algebra of functions vanishing at infinity on the fibres. This is given by taking $f_0\in C_0(\mathbb{R},A)$.

\noindent b)\quad The unital Thom algebra, given by adjoining $\mathbb{C}$ to the first case (a)
(i.e.\  the one point compactification).

\noindent c)\quad The unital associated circle bundle algebra.  This is given by taking $f_0\in C(\mathbb{R}_\infty,A)$. 
\end{theorem}
\noindent {\bf Proof:}\quad Use Propositions~\ref{vcahjgkvyfdf1} and \ref{vcahjgkvyfdf2}, and take completions. \quad$\square$

\medskip Note that the whole idea of a star operation on $T_{\mathbb{N}}(L)$ depends on the existence of a star operaton $\star:L\to \overline{L}$.
So when do we have a star operation $\star:L\to \overline{L}$? Note that for $L$ to be isomorphic to $\overline{L}$ as a bimodule is likely not sufficient, we need to have $\mathrm{bb}=\overline{\star}\,\star:L\to \overline{\overline{L}}$. But let us simply suppose that $L$ and $\overline{L}$ are isomorphic, and see how far we can get.

\begin{propos}\label{cvgadfsfg}
Suppose that $L$ is a left line bimodule and that $\theta:L\to \overline{L}$ is an invertible bimodule map. 
Then the map 
\begin{eqnarray} \label{vcakvkchc1}
L\stackrel{\theta}\longrightarrow \overline{L}
\stackrel{\overline{\theta}}\longrightarrow \overline{\overline{L}}
\stackrel{\mathrm{bb}^{-1}}\longrightarrow L
\end{eqnarray}
is given by $e\mapsto e.z$ where $\Phi_L(z^*)=z\in Z(A)$ is invertible. A bimodule map $\star:L\to\overline{L}$ satisfying $\overline{\star}\,\star=\mathrm{bb}:L\to  \overline{\overline{L}}$ exists if and only if there is a $y\in Z(A)$ with $\Phi_L(y^*)\,y=z^{-1}$, in which case we can have $\star(e)=\theta(e).y$. 
\end{propos}
{\bf Proof:}\quad 
Write $\theta(x)=\overline{g(x)}$. Then the composition in (\ref{vcakvkchc1})
is a bimodule map from $L$ to itself, and by Proposition~\ref{cagjcjxdz} is given by $e\mapsto e.z$,
for some $z\in Z(A)$. Consider the following map:
\begin{eqnarray}  \label{vcakvkchc}
L\stackrel{\theta}\longrightarrow \overline{L}
\stackrel{\overline{\theta}}\longrightarrow \overline{\overline{L}}
\stackrel{\overline{\overline{\theta}}}\longrightarrow \overline{\overline{\overline{L}}}
\stackrel{\overline{\mathrm{bb}^{-1}}}\longrightarrow \overline{L}
\stackrel{\theta^{-1}}\longrightarrow L\ .
\end{eqnarray}
One way of calculating (\ref{vcakvkchc}) is as
\begin{eqnarray*}
L\stackrel{\theta}\longrightarrow \overline{L}
\stackrel{\overline{\mathrm{bb}^{-1}\,\overline{\theta}\,\theta}}\longrightarrow \overline{L}
\stackrel{\theta^{-1}}\longrightarrow L\ ,
\end{eqnarray*}
which gives $e\mapsto\overline{e'}\mapsto\overline{e'.z}=z^*\,\overline{e'}\mapsto
z^*\,e$. Alternatively the composition of the last three maps in (\ref{vcakvkchc}) is 
\begin{eqnarray*}
\overline{\overline{e}} \longmapsto \overline{\overline{\overline{g(e)}}}
\longmapsto \overline{g(e)}\longmapsto e\ ,
\end{eqnarray*}
so (\ref{vcakvkchc}) is just another way of writing (\ref{vcakvkchc1}). 
Comparing these gives $z^*\,e=e\,z$ for all $e\in L$, so $\Phi_L(z^*)=z$. Any map 
$\star:L\to\overline{L}$ is given by $\star(e)=\theta(e).y$, for some $y\in Z(A)$, hence
\begin{eqnarray*}
\overline{\star}\,\star(e) &=& \overline{\star(g(e))}.y \cr
&=& (\overline{\overline{g(g(e))}.y}).y \cr
&=& y^*.\overline{\overline{g(g(e))}}.y\ ,
\end{eqnarray*}
and then
\begin{eqnarray*}
\mathrm{bb}^{-1}\,\overline{\star}\,\star(e) &=& y^*\,g(g(e))\,y\ =\ 
g(g(e))\,\Phi_L(y^*)\,y\cr
&=& e\,z\,\Phi_L(y^*)\,y\ .\quad\square
\end{eqnarray*}

\medskip So just how restrictive is the condition in Proposition~\ref{cvgadfsfg} for the existence of a star structure on a line module $L$, given that $L$ and $\overline{L}$ are in the same bimodule isomorphism class? Suppose that $A$ is a unital $C^*$-algebra, and then its centre $Z(A)$ 
consists of complex valued continuous functions on a compact topological space. If we restrict to the case where $\Phi_L$ is the identity, then the $z$ of Proposition~\ref{cvgadfsfg} is a Hermitian function, i.e.\ a nowhere vanishing continuous real valued function on $X$. Hence the existence of $y$ reduces to asking whether $z$ is positive or not, looking at components of $X$.

\section{An example of constructing line modules}\label{constcxagjy}
Here we take an example of a construction of line modules, which has been used in the literature
\cite{LandiDecon,Landi4DHopf}. Suppose that $B$ is a comodule algebra for the coaction
$\rho:B\to B\tens \mathbb{C}G$ of the group algebra of a (discrete) group $G$. Let $A$ be the invariant part of $B$ under the coaction. Further suppose that there are two column vectors
\begin{eqnarray}
\underline{v} \ =\ \left(\begin{array}{c}v_1 \\v_2 \\\vdots \\v_n\end{array}\right)\ ,\quad
\underline{w} \ =\ \left(\begin{array}{c}w_1 \\w_2 \\\vdots \\w_n\end{array}\right)\ ,
\end{eqnarray}
with the following properties. All the $v_i$ are in a particular $g\in G$ graded part $B_g$
of $B$ (i.e.\ $\rho(v_i)=v_i\tens \underline{g}$ where $\underline{g}$ is a basis element of the group algebra $\mathbb{C}G$). All the $w_i$ are in the $g^{-1}\in G$ graded part of $B$. Finally, the column vectors satisfy the following equation
\begin{eqnarray*}
\underline{w}^T \,\underline{v} \ =\ \sum_i w_i\,v_i\ =\ 1_B\ .
\end{eqnarray*}
Then $P=\underline{v}\,\underline{w}^T$ is an $n$ by $n$ matrix with entries in $A$
and satisfies $P^2=P$. 
In the case where $g\in G$ is the group identity, the corresponding line module would be trivial as a module over $A$, but not in general. This construction was used in  
\cite{LandiDecon} to explain the classical Dirac monopole bundle, and it was subsequently used to construct very non-trivial noncommutative bundles for example in \cite{ncgaugeBM} or \cite{Landi4DHopf}. 

Let $L$ and $L^\circ$ be defined as
\begin{eqnarray}
L &=& \big\{ b.\underline{w}^T \subset A^{\oplus n} \ \big|\  b\in B_g\big\}\ ,\cr
L^\circ &=& \big\{ \underline{v}.b \subset A^{\oplus n} \ \big|\  b\in B_{g^{-1}}\big\}\ .
\end{eqnarray}
These are $A$ bimodules, with left and right action
\begin{eqnarray}
&&a\,\la\,(b.\underline{w}^T) \ =\ ab.\underline{w}^T\ ,\quad
( b.\underline{w}^T) \,\ra\, a\ =\ b.\underline{w}^T.P(a)\ =\ ba.\underline{w}^T\ ,\cr
&&a\,\la\,( \underline{v}.b) \ =\ P(a). \underline{v}.b\ =\ \underline{v}.ab\ ,\quad
( \underline{v}.b) \,\ra\, a\ =\ \underline{v}.ba\ ,
\end{eqnarray}
where $P(a)=\underline{v}\,a\,\underline{w}^T$. The evaluation map can be quite simply defined by matrix multiplication:
\begin{eqnarray}
\ev:L\tens_A L^\circ\to A\ ,\quad \ev(b.\underline{w}^T\tens \underline{v}.b')\ =\ 
b.\underline{w}^T\, \underline{v}.b'\ =\ b\,b'\ .
\end{eqnarray}
For the coevaluation map, we choose $c\tens c'\in B_{g^{-1}}\tens B_g$ (summation implicit) so that $c\,c'=1$. Of course, we can choose $c\tens c'=w_i\tens v_i$, but it may cause less complications to keep $c\tens c'$ separate, as they will frequently appear in the same formula as $\underline{v}$ and $\underline{w}$.  Then define
\begin{eqnarray}
\coev:A\to L^\circ\tens_A L\ ,\quad \coev(a)\ =\ \underline{v}\,a\,c\tens c'\,\underline{w}^T\ .
\end{eqnarray}
Now we check the required properties for the evaluation and coevaluation:
\begin{eqnarray*}
(\ev\tens\id_L)(b.\underline{w}^T\tens\coev(1)) &=& b.\underline{w}^T.\underline{v}\,c\,\la\, c'\,\underline{w}^T \cr
&=& b\,c\,c'\,\underline{w}^T\ =\ b.\underline{w}^T\ ,\cr
(\id_{L^\circ}\tens\ev)(\coev(1)\tens  \underline{v}\,b ) &=&  \underline{v}\,c\,\ra\, c'\,\underline{w}^T\,\underline{v}.b \cr
&=& \underline{v}\,c\, c'\,b\ =\  \underline{v}\,b\ .
\end{eqnarray*}

\begin{propos}
Suppose that the $G$-graded algebra $B$ is a $\mathbb{C} G$ Hopf-Galois extension of $A$.Then $L$
is a left line module.
\end{propos}
\noindent {\bf Proof:}\quad As $B$ is a $\mathbb{C} G$ Hopf-Galois extension of $A$, the product
$B_g\tens B_{g^{-1}}\to B_0=A$ is surjective. Then we can choose $c''\tens c'''\in B_g\tens B_{g^{-1}}$ (summation implicit) with $c''\,c'''=1$.
We shall define maps which we will then check are inverses to evaluation and coevaluation, by
\begin{eqnarray*}
\coev^{-1}(\underline{x}\tens\underline{y}^T) &=& \underline{w}^T\,\underline{x}\,\underline{y}^T\,\underline{v}\ ,\cr
\ev^{-1}(a) &=& a\,c''\underline{w}^T\tens \underline{v}\,c'''\ .
\end{eqnarray*}
First we do the easy checks:
\begin{eqnarray*}
\coev^{-1}\circ\coev(1) &=& \coev^{-1}(\underline{v}\,c\tens c'\,\underline{w}^T) \cr
&=&  \underline{w}^T\, \underline{v}\,c\, c'\,\underline{w}^T\,\underline{v} \ =\ 1\ ,\cr
\ev\circ\ev^{-1}(a) &=& \ev(a\,c''\underline{w}^T\tens \underline{v}\,c''')\ =\ a\,c''.c'''\ =\ a\ .
\end{eqnarray*}
The more difficult checks depend on the Hopf-Galois condition:
\begin{eqnarray*}
\coev\circ\coev^{-1}(\underline{v}\,b\tens b'\,\underline{w}^T) &=&\coev( \underline{w}^T\,\underline{v}\,b\, b'\,\underline{w}^T\,\underline{v}) \cr
&=& \underline{v}\, \underline{w}^T\,\underline{v}\,b\, b'\,\underline{w}^T\,\underline{v}\,c\tens c'\,\underline{w}^T \cr
&=& \underline{v}\, b\, b'\,c\tens c'\,\underline{w}^T
\end{eqnarray*}
At this point we use the fact that product gives an isomorphism $B_{g^{-1}}\tens_A B_g\to A$, as this shows that $b\, b'\,c\tens_A c'=b\tens_A b'$. Next we have
\begin{eqnarray*}
\ev^{-1}\circ\ev(b'.\underline{w}^T\tens \underline{v}.b) &=& \ev^{-1}(b'\,b) \cr
&=& b'\,b\, c''\underline{w}^T\tens \underline{v}\,c''' \ .
\end{eqnarray*}
Now we use the fact that product gives an isomorphism $B_{g}\tens_A B_{g^{-1}}\to A$, as this shows that $b'\,b\, c''\tens_A c'''=b'\tens_A b$. Finally we write $b'\,b\, c''\tens c'''-b'\tens b=
r.a\tens s-r\tens a.s$ for $a\in A$ (summation implicit), and then in $L\tens_A L^\circ$,
\begin{eqnarray*}
&& \ev^{-1}\circ\ev(b'.\underline{w}^T\tens \underline{v}.b) - b'.\underline{w}^T\tens \underline{v}.b\cr
&=& r\,a\,\underline{w}^T\tens \underline{v}\,s -r\,\underline{w}^T\tens \underline{v}\,a\,s \cr
&=& r\,\underline{w}^T\,P(a)\tens \underline{v}\,s -r\,\underline{w}^T\tens P(a)\,\underline{v}\,s\ =\ 0 \ .\quad\square
\end{eqnarray*}

\medskip
For a coaction of a Hopf $*$-algebra on a comodule $V$, the coaction on the conjugate comodule is defined by $\overline{e}\mapsto \overline{e_{[0]}}\tens e_{[1]}^*$. The usual star algebra structure on the group algebra of $G$ is such that every group element is unitary, i.e.\ $g^*=g^{-1}$ for $g\in G$.

\begin{propos}
Suppose that the $G$-graded algebra $B$ is a star algebra with $b\in B_g$ implying
$b^*\in B_{g^{-1}}$. 
If $\underline{w}^T=\underline{v}^*$ (star on a matrix being star element-wise, then transpose), then we can define a non-degenerate Hermitian inner product by
$G(\overline{  \underline{y}^T  }) = (\underline{y}^T)^*$. 
In addition the corresponding inner product is positive, meaning that each $\<\underline{y}^T,\overline{\underline{y}^T}\>$ is a sum of elements of the form $a^*\,a$ for $a\in A$. 
\end{propos}
\noindent {\bf Proof:}\quad First we check that $G:\overline{L}\to L^\circ$, by
\begin{eqnarray*}
G(\overline{b\,\underline{w}^T})\ =\ (b\,\underline{w}^T)^*\ =\ \underline{v}^{**}\,b^*\ =\ \underline{v}\,b^*\ \in L^\circ\ .
\end{eqnarray*}
Now we check that it is a bimodule map. For $a\in A$,
\begin{eqnarray*}
G(a\,\la\,\overline{b\,\underline{w}^T}) &=& G(\overline{(b\,\underline{w}^T)\,\ra \, a^*}) 
\ =\ G(\overline{b\,a^*\,\underline{w}^T})
\cr &=& \underline{v}\,(b\,a^*)^*\ =\ \underline{v}\,a\,b^*\ =\ a\,\la\,(\underline{v}\,b^*)\ , \cr                 
G(\overline{b\,\underline{w}^T}\,\ra\,a) &=& G(\overline{a^*\,\la\,(b\,\underline{w}^T)})\ =\
G(\overline{a^*\,b\,\underline{w}^T})
\cr &=& \underline{v}\, (a^*\,b)^*\ =\  \underline{v}\, b^*\,a\ =\ (\underline{v}\, b^*)\,\ra\,a\ .
\end{eqnarray*}
The symmetry of the corresponding inner product is checked as follows:
\begin{eqnarray*}
\<\underline{y}^T,\overline{\underline{x}^T}\>^* &=& \ev(\underline{y}^T\tens (\underline{x}^T)^*)^*\cr
&=& (\underline{y}^T\, (\underline{x}^T)^*)^* \cr
&=& \Big(\sum y_i\,x_i^*\Big)^*\ =\ \sum x_i\,y_i^*\cr
&=& \<\underline{x}^T,\overline{\underline{y}^T}\>\ .
\end{eqnarray*}
The formula for $G^{-1}:L^\circ\to \overline{L}$ is just $G^{-1}(\underline{z})=\overline{\underline{z}^*}$. To check positivity, we use the formula above
\begin{eqnarray*}
\<\underline{y}^T,\overline{\underline{y}^T}\>\ =\ \sum y_i\,y_i^*\ .\quad\square
\end{eqnarray*}

\medskip In the beginning of this section, using the vectors $\underline{v}$ and
$\underline{w}$ was sold as a method of constructing examples of line modules. In fact,
it is rather more than that, any line module can be constructed by this method. The catch is that different bundles may require different Hopf-Galois extensions to realise them. An interesting question would be whether there is some form of universal Hopf-Galois extension from which any line module can be constructed up to isomorphism -- see the comment in Section~\ref{pichgcjf}.

\begin{propos}
Given any left line module $L$ over an algebra $A$, there is an integer graded Hopf-Galois extension $C$ of $A$ so that $L$ is given, up to isomorphism, by the vector construction in the beginning of this section.
\end{propos}
\noindent {\bf Proof:}\quad From Theorem~\ref{vchuaytdrtys}, there is an integer graded Hopf-Galois extension $C$ of $A$ so that, by construction, $C_1=L$, $C_{-1}=L^\circ$
and $C_0=A$. As multiplication  $:C_{1}\tens C_{-1}\to C_0$ is onto, there are elements $w_i\in C_{-1}$ and $v_i\in C_1$ ($1\le i\le n$) so that
\begin{eqnarray*}
\sum w_i\,v_i\ =\ 1\ .
\end{eqnarray*}
The bimodules $C_1\,\underline{w}^T$ and $L$ are isomorphic
by the maps
\begin{eqnarray*}
b\in C_1\mapsto b\,\underline{w}^T\ ,\quad \underline{y}^T\mapsto
\underline{y}^T\,\underline{v}\ .\quad\square
\end{eqnarray*}

\medskip It is worthwhile to note that the vectors $\underline{v}$ and $\underline{w}$ do not alter the isomorphism class of the line module, that is determined purely by the graded algebra. However their existence is used to demonstrate that we actually have a line module, and in constructing the projection matrix
$P=\underline{v}\,\underline{w}^T$.

\begin{example} \label{vfauxkyytuyt}
The example of a line bundle on $\mathbb{C}_q[SL_2]$ is well known; see \cite{HajMaj:pro}. Suppose that $q\in\mathbb{C}$ with $q^2\ne 1$. 
The quantum group $\C_q[SL_2]$ has generators $a,b,c,d$
 with relations:
\[ ba=qab\ ,\ ca=qac\ ,\ db=qbd\ ,\ dc=qcd\ , \ cb=bc \ , \]
\[
 da-ad=q(1-q^{-2}) bc \ , \ ad-q^{-1}bc=1.\] 
The coproduct $\Delta$ and counit $\eps$ have the usual matrix
coalgebra form. We denote the antipode or `matrix inverse' by $S$:
\begin{eqnarray*}
S\left(\begin{array}{cc}a & b \\c & d\end{array}\right) \,=\,
\left(\begin{array}{cc}d & -q\,b \\ -q^{-1}\,c & a\end{array}\right)\ .
\end{eqnarray*} 
The algebra $\C_q[SL_2]$ equipped with the star operation $a^*=d$, $d^*=a$, $c^*=-q\,b$
and $b^*=-q^{-1}c$, where $q$ is real, is denoted $\C_q[SU_2]$.

There is a grading on $\mathbb{C}_q[SL_2]$ for which the generators $a,c$ have degree $+1$
and $b,d$ have degree $-1$. Now we set
\begin{eqnarray*}
\underline{w}\ =\ \left(\begin{array}{c}a \\c/q\end{array}\right)\ ,\quad
\underline{v}\ =\ \left(\begin{array}{c}d \\-b\end{array}\right)\ ,
\end{eqnarray*}
and then $(\underline{w}^T)^*=\underline{v}$ and $\underline{w}^T\,\underline{v}=1$. It follows that we have Hermitian metrics. 

But what of a star operation $\star:L\to \overline{L}$? We do not expect one in the case of the sphere. To explain, we can consider constructing line bundles on the ordinary sphere. The complex line bundles are given by `clutching functions' from the equator to $\mathbb{C^*}$ (or $S^1$ with a metric) (see \cite{AtiK}). These are classified by the winding number in $\mathbb{Z}$. If $L$ has winding number $n$, then 
$L^\circ$ has winding number $-n$, so we do not expect an isomorphism from $L$ to $\overline{L}$. Looked at another way, real line bundles (with metric) are classified by clutching functions from the equator to $\mathbb{Z}/2$, giving only one real bundle (the trivial one). 
\end{example}

\section{An example of a Chern class in de Rham cohomology}\label{acjfjdssd}
Here we shall consider the differential geometry of a line module, and do several calculations.
It has not escaped our notice that connections on line bundles is an important part of gauge theory in theoretical physics, and the framework of metrics which we have presented, and the space of covariant derivatives in this section, may be relevant (see \cite{ncgaugeBM}). 
In particular we shall be concerned with the definition in \cite{KobNom} of Chern class of a vector bundle as a de Rham cohomology class given by taking traces of powers of the curvature. 
 In particular, in \cite{KobNom} there is a proof that the de Rham cohomology class does not depend on the covariant derivative $\nabla$ chosen on the vector bundle. First this depends on the fact that the set of covariant derivatives is connected (in fact, it is an affine space), and then takes differentiable paths between different covariant derivatives (parameterised by time $t$, say). The derivative of the constructed element of $\Omega^2$ is shown to be in the image of $\extd:\Omega^1\to\Omega^2$, so the de Rham class for the two covariant derivatives is the same. We shall attempt to produce a `plausible' copy of this proof in a noncommutative context, for an example with line modules. 

Suppose that we have a differential calculus $(\Omega^*A,\extd)$ on the algebra $A$. In this section we want to study the left covariant derivatives on the line module $L=B_g$. 
As in Section~\ref{constcxagjy}, suppose that $B$ is a Hopf-Galois extension of the algebra $A$ for the coaction
$\rho:B\to B\tens \mathbb{C}G$ of the group algebra of a (discrete) group $G$. It will be convenient to fix
an element $g\in G$,   
$c\tens c'\in 
B_{g^{-1}}\tens_{\mathbb{C}} B_g$ (summation implicit) so that $c\,c'=1_B$ and $c''\tens c'''\in 
B_{g}\tens_{\mathbb{C}} B_{g^{-1}}$ (summation implicit) so that $c''\,c'''=1_B$. In the star algebra case we shall also assume, without further loss of generality, that $c\tens c'={c'}^*\tens c^*\in B_{g^{-1}}\tens_{\mathbb{C}} B_g$. This can be done by taking a new $c\tens c'$ to be the average of the old
$c\tens c'$ and ${c'}^*\tens c^*$. We also use $h\tens h'$ and $f\tens f'$ as independent copies of $c\tens c'$.

\begin{propos}  \label{vchjakjxcjj}
A general left covariant derivative on $L=B_g$ can be written as
\begin{eqnarray*} 
\nabla_L(e) &=& \extd( e\,c)\tens c' +
\tilde\Gamma(e)\ \in\ \Omega^1 A\tens_A L\ ,
\end{eqnarray*}
where $\tilde\Gamma:L\to \Omega^1 A\tens_A L$ is a left $A$-module map. 
\end{propos}
\noindent {\bf Proof:}\quad If we define $\tilde\Gamma(e)=\nabla_L(e) - \extd( e\,c)\tens c'$, then for $a\in A$,
\begin{eqnarray*}
\tilde\Gamma(a.e)  &=& \nabla_L(a.e) - \extd( a\,e\,c)\tens c'  \cr
 &=& \extd a\tens e +a.\nabla_L(e) - \extd a.e\,c\tens c' -  a.\extd(e\,c)\tens c'  \cr
 &=& \extd a\tens e - \extd a.e\,c\tens c' + a.\tilde\Gamma(e) \ .
\end{eqnarray*}
Now we use the fact that for the Hopf-Galois extension, $e\,c\tens c'=1_A\tens e\in A\tens_A B_g$.\quad$\square$

\medskip Suppose that $B$ is equipped with a differential calculus which is also $G$-graded, and that 
$B_g.\Omega^1 A.B_{g^{-1}}\subset \Omega^1 A$. This is not an arbitrary condition, it fits with the idea that the 1-forms on $A$ are the horizontal invariant 1-forms on $B$, see \cite{ncsheaf}. Using this we can define a useful bimodule map
\begin{eqnarray} \label{zetmghjchdt}
\sigma_0:L\tens_A\Omega^1 A\to \Omega^1 A\tens_A L\ ,\quad e\tens\eta\mapsto e.\eta.c\tens c'\ .
\end{eqnarray}
This is obviously a left module map, and to show that it is a right module map
 we use $a\,c\tens c'=c\tens c'\,a\in B_{g^{-1}}\tens_A B_g$ for all $a\in A$.
 Also $\sigma_0$ is invertible, with inverse $\sigma_0^{-1}(\eta\tens b)=c''\tens c'''.\eta.b$.
  This will prove useful in the next result:

\begin{propos} \label{cvajkhxsj}
The general form of $\tilde{\Gamma}$ is given by the following formula, for some $\zeta\in \Omega^1 A$,
\begin{eqnarray*}
\tilde{\Gamma}(b) \ =\ e.\zeta.c\tens c'   \ ,
\end{eqnarray*}
and this gives the general form of the covariant derivative and curvature as
\begin{eqnarray*}
\nabla_L(e) &=& \extd(e\,c)\tens c' + e.\zeta.c\tens c'\cr
\mathrm{R}_L(e) &=& e\,h.\extd(h'.\zeta.c)\tens c'- e\,f.\extd(f'\,c)\wedge\extd(c'\,h)\tens h'    - e.\zeta.c\wedge c'.\zeta.h\tens h'   \cr
&=&  e\,h.\extd(h'.\zeta.c)\tens c'-\extd(e\,c)\wedge\extd(c'\,h)\tens h'    - e.\zeta.c\wedge c'.\zeta.h\tens h'\ . 
\end{eqnarray*}
\end{propos}
\noindent {\bf Proof:}\quad 
 To establish the form of $\tilde{\Gamma}$
we use the invertible bimodule map $\sigma_0$ from (\ref{zetmghjchdt}). 
In view of Proposition~\ref{xdcxyjtalj} (for $F=\Omega^1 A$) we can write the left module map
\begin{eqnarray*}
\sigma_0^{-1}\,\tilde{\Gamma}:L\to L\tens_A\Omega^1 A
\quad\mathrm{as}\quad e\mapsto e\tens\zeta\ ,
\end{eqnarray*}
 for some $\zeta\in \Omega^1 A$, and thus $\tilde{\Gamma}$ has the required form. The curvature is calculated as
\begin{eqnarray*}
\mathrm{R}_L(e) &=&
(\extd\tens\id-\id\wedge\nabla_L)\,\nabla_L(e) \cr
&=& (\extd\tens\id-\id\wedge\nabla_L)\,
( \extd(e\,c)\tens c' + e.\zeta.c\tens c' )\cr
&=& \extd(e.\zeta.c)\tens c'-\extd(e\,c)\wedge\nabla_L( c')  - e.\zeta.c\wedge\nabla_L( c') \cr
&=& \extd(e.\zeta.c)\tens c'-\extd(e\,c)\wedge\extd(c'\,h)\tens h'-\extd(e\,c)\wedge c'.\zeta.h\tens h'\cr
&& - \ e.\zeta.c\wedge  \extd(c'\,h)\tens h'     - e.\zeta.c\wedge c'.\zeta.h\tens h'\ .
\end{eqnarray*}
To simplify this, we first observe that in $\Omega^* A \tens_{\mathbb{C}} L$,
\begin{eqnarray*}
e.\zeta.c\tens c'  &=& e.\zeta.c\,c'\,h\tens h'    \ ,\cr
\extd(e.\zeta.c)\tens c' &=&   \extd(e.\zeta.c).c'\,h\tens h' - e.\zeta.c\wedge\extd(c'\,h)\tens h'  \ ,
\end{eqnarray*}
and on taking this in $\Omega^2 A \tens_A L$ we see that $e.\zeta.c\wedge\extd(c'\,h)\tens h'=0$. 
Next in $\Omega^* A \tens_{\mathbb{C}} L$,
\begin{eqnarray*}
e.\zeta.c\tens c' &=& e\,h\,h'.\zeta.c\tens c' \ ,\cr
\extd(e.\zeta.c)\tens c' &=& \extd(e\,h)\wedge h'.\zeta.c\tens c' + e\,h.\extd(h'.\zeta.c)\tens c' \ ,
\end{eqnarray*}
and using these results we can rewrite the curvature as
\begin{eqnarray*}
\mathrm{R}_L(e) 
&=& e\,h.\extd(h'.\zeta.c)\tens c'-\extd(e\,c)\wedge\extd(c'\,h)\tens h'    - e.\zeta.c\wedge c'.\zeta.h\tens h'\ . 
\end{eqnarray*}
Finally we compute
\begin{eqnarray*}
&& \extd(e\,c)\wedge\extd(c'\,h)\tens h' \cr
&=& \extd(e\,f\,f'\,c)\wedge\extd(c'\,h)\tens h' \cr
&=& \extd(e\,f).f'\,c\wedge\extd(c'\,h)\tens h' + e\,f.\extd(f'\,c)\wedge\extd(c'\,h)\tens h' \cr
&=& \extd(e\,f)\wedge f'\,c.\extd(c'\,h)\tens h' + e\,f.\extd(f'\,c)\wedge\extd(c'\,h)\tens h' \cr
&=& \extd(e\,f)\wedge \extd(f'\,c\,c'\,h)\tens h' -\extd(e\,f)\wedge\extd(f'\,c).c'\,h\tens h' + e\,f.\extd(f'\,c)\wedge\extd(c'\,h)\tens h' \cr
&=& \extd(e\,f)\wedge \extd(f'\,h)\tens h' -\extd(e\,f)\wedge\extd(f'\,c)\tens c'\,h\,h' + e\,f.\extd(f'\,c)\wedge\extd(c'\,h)\tens h' \cr
&=&  e\,f.\extd(f'\,c)\wedge\extd(c'\,h)\tens h' \ .\quad\square
\end{eqnarray*}

\medskip
In \cite{BegMa4} the condition for a left covariant derivative to preserve a Hermitian inner product is given, under the assumption that the covariant derivative is actually a bimodule covariant derivative. 
It turns out that the condition can be stated in a form which does not make use of this bimodule covariant derivative assumption. The reason is that a left covariant derivative $\nabla_L$ on a module $E$ automatically gives a right covariant derivative $\check\nabla$ on the conjugate module $\overline{E}$, so we can form a covariant derivative $\overline{E}\tens_A E\to \overline{E}\tens_A \Omega^1 A \tens_A E$. This is effectively what we do in Proposition~\ref{ycfuxsarr}, and is the reason why we do not actually require a bimodule covariant derivative.

\begin{propos} \label{ycfuxsarr}
Suppose that we define a metric using the star operation $\star:B_g\to \overline{B_{g^{-1}}}$
as $G(\overline{b})=b^*$. 
The equation for the covariant derivative in Proposition~\ref{cvajkhxsj} to preserve the metric is that $\zeta+\zeta^*=0$.
\end{propos}
\noindent {\bf Proof:}\quad We need to solve the following equation,
\begin{eqnarray}   \label{cvgajsasrty1}
0 &=& (\check\nabla\tens\id+\id\tens\nabla_{L})\,(G^{-1}\tens\id_L)\coev_L(1_A) \cr
&=& (\check\nabla\tens\id+\id\tens\nabla_{L})\,(\overline{c^*}\tens c')\cr
&=& \check\nabla(\overline{c^*})\tens c' +
\overline{c^*}\tens \nabla_L(c')\ .
\end{eqnarray}
We use the definition for $\check\nabla$ given in \cite[Proposition~3.2]{BegMa4},
\begin{eqnarray*}
\check\nabla (\overline{b}) &=& (\id\tens\star^{-1})
\Upsilon\,\overline{\nabla_L(b)}\cr 
&=&  (\id\tens\star^{-1})  \Upsilon\, \big( \overline{     \extd (b\,h)\tens h'+ b.\zeta.h\tens h'
}\big)  \cr
&=& \overline{h'}\tens \extd(h^*\, b^*)+ \overline{h'}\tens h^*.\zeta^*.b^*\ .
\end{eqnarray*}
If we substitute this into (\ref{cvgajsasrty1}), we get 
\begin{eqnarray*}
0 &=&   \overline{h'}\tens \extd(h^*\, c)\tens c'+ \overline{h'}\tens h^*.\zeta^*.c\tens c'\cr
&&+\ \overline{c^*}\tens \extd (c'\,h)\tens h'
+ \overline{c^*}\tens c'.\zeta.h\tens h'\ .
\end{eqnarray*}
Applying $\star^{-1}\tens\id\tens\id$ and relabeling gives
\begin{eqnarray*}
0 &=&   {c'}^*\tens \extd(c^*\, h)\tens h'+ {c'}^*\tens c^*.\zeta^*.h\tens h'\cr
&&+\ c\tens \extd (c'\,h)\tens h'
+ c\tens c'.\zeta.h\tens h'\ .
\end{eqnarray*}
Next we use the assumption that $c\tens c'={c'}^*\tens c^*$ to write this as
\begin{eqnarray*}
0 &=& c\tens (2\,\extd(c'\,h)+c'.(\zeta+\zeta^*).h   )\tens h'\ .
\end{eqnarray*}
Note that there is an isomorphism from $B_{g^{-1}}\tens_A \Omega^1 A \tens_A B_g$ to 
$\Omega^1 A$ given by taking the product of all the factors. The inverse is just $\eta\in \Omega^1 A$ mapping to $c\tens c'.\eta.h\tens h'$. Using this isomorphism we can restate our condition as
\begin{eqnarray*}
0 &=& 2\,c.\extd(c'\,h).h'+\zeta+\zeta^*\ \in \ \Omega^1 A  \ .
\end{eqnarray*}
Finally, 
\begin{eqnarray*}
c.\extd(c'\,h).h' &=& f\,f'\,c.\extd(c'\,h).h' \cr
&=& f.\extd(f'\,c\,c'\,h).h'-f.\extd(f'\,c).c'\,h\,h'  \cr
&=& f.\extd(f'\,h).h'-f.\extd(f'\,c).c' \ =\ 0  \ .\quad\square
\end{eqnarray*}

\medskip We can suppose that $B$ has a differential calculus $(\Omega^*B,\extd,\wedge)$
which contains $(\Omega^*A,\extd,\wedge)$ as a sub-differential graded algebra. We shall assume this in the rest of this section, and it will be true in Example~\ref{ctfuxkklyfsr}. In particular we have a useful 1-form $\kappa\in\Omega^1 B$, defined by 
\begin{eqnarray}\label{vckuacsxs}
\kappa\ =\ c.\extd c'\ .
\end{eqnarray}
Any $a\in A$ commutes with $\kappa$ up to an element of $\Omega^1 A$, as can be seen from the formula
\begin{eqnarray}
c.\extd (c'\,a\,h).h'\ =\ \kappa.a-a.\kappa+\extd a\ .
\end{eqnarray}

\begin{cor}\label{dxfghjjmxfhj}
Using the formula for the covariant derivative in Proposition~\ref{cvajkhxsj}, and supposing that $B$ has a differential calculus containing the differential calculus for $A$ as a subalgebra, we have the following formulae for $\sigma$ (if it exists, see Definition~\ref{ppll}):
\begin{eqnarray*}
\sigma(e\tens\extd a) &=& e.\extd(a).c\tens c' + e\,(a.(\zeta-\kappa)-(\zeta-\kappa).a)\,c\tens c'\ .
\end{eqnarray*}
\end{cor}
\noindent {\bf Proof:}\quad Just using $\Omega^1 A$,
\begin{eqnarray*}
\nabla_L(e\,a)-\nabla_L(e).a &=& \extd(e\,a\,c)\tens c'-\extd(e\,c)\tens c'\,a+
e\,a.\zeta.c\tens c'  -e.\zeta.c\tens c'\,a   \cr
&=& \extd(e\,h\,h'\,a\,c)\tens c'-\extd(e\,c)\tens c'\,a +e\,(a.\zeta-\zeta.a)\,c\tens c' \cr
&=& e\,h.\extd(h'\,a\,c)\tens c'+ \extd(e\,h).h'\,a\,c\tens c' -\extd(e\,c)\tens c'\,a \cr
&& +\ e\,(a.\zeta-\zeta.a)\,c\tens c' \cr
&=& e\,h.\extd(h'\,a\,c)\tens c'+ e\,(a.\zeta-\zeta.a)\,c\tens c' \ .
\end{eqnarray*}
Now we continue using the differential calculus on $B$,
\begin{eqnarray*}
\nabla_L(e\,a)-\nabla_L(e).a
&=& e\,h.\extd(h').a\,c\tens c'+e\,h\,h'.\extd(a).c\tens c'  \cr
&& +\ e\,h\,h'\,a.\extd(c)\tens c'+ e\,(a.\zeta-\zeta.a)\,c\tens c'  \cr
&=& e\,h.\extd(h').a\,c\tens c'+e.\extd(a).c\tens c'  \cr
&& +\ e\,a.\extd(c).c' \,h\tens h'+ e\,(a.\zeta-\zeta.a)\,c\tens c'  \cr
&=& e\,h.\extd(h').a\,c\tens c'+e.\extd(a).c\tens c'  \cr
&& -\ e\,a.h.\extd(h'). c\tens c'+ e\,(a.\zeta-\zeta.a)\,c\tens c' \ . \quad\square
\end{eqnarray*}

\begin{cor}\label{dxfghxfhj}
Using the formula for the covariant derivative in Proposition~\ref{cvajkhxsj}, and supposing that $B$ has a differential calculus containing the differential calculus for $A$ as a subalgebra, we have the following formulae for the curvature $\mathrm{R}_L$:
\begin{eqnarray*}
\mathrm{R}_L(e) &=& e.\big( 
\extd(\zeta-\kappa)-(\zeta-\kappa)\wedge(\zeta-\kappa) \big).c\tens c'\ .
\end{eqnarray*}
\end{cor}
\noindent {\bf Proof:}\quad From Proposition~\ref{cvajkhxsj},
\begin{eqnarray*}
\mathrm{R}_L(e) &=& e\,h.\extd(h'.\zeta.c)\tens c'-\extd(e\,c)\wedge\extd(c'\,h)\tens h'    - e.\zeta.c\wedge c'.\zeta.h\tens h'   \cr
&=&  e\,\kappa\wedge\zeta.c\tens c'+e.\extd\zeta.c\tens c'-e.\zeta\wedge\extd c\tens c'\cr
&& -\ \extd(e\,c)\wedge\extd(c'\,h)\tens h'    - e.\zeta\wedge \zeta.h\tens h' \ .
\end{eqnarray*}
We use
\begin{eqnarray*}
  -\, \extd(e\,c)\wedge\extd(c'\,h)\tens h'    
&=&  -\, \extd e.c\wedge\extd c'.h\tens h'  - \extd e.c\wedge c'.\extd h\tens h' \cr
&&  -\, e.\extd c\wedge\extd c'.h\tens h'  - e.\extd c\wedge c'.\extd h\tens h' \cr
&=&  -\, \extd e\wedge c.\extd c'.h\tens h'  - \extd e.\wedge \extd h\tens h'\,c\,c' \cr
&&  -\, e.\extd c\wedge\extd c'.h\tens h'  - e.\extd c.c'\wedge \extd h.h'\,f\tens f' \cr
&=& -\, e.(\extd\kappa+\kappa\wedge\kappa).c\tens c' \ .\quad\square
\end{eqnarray*}

\medskip We define the trace $\omega_L\in\Omega^2 A$ of the curvature $\mathrm{R}_L:L\to \Omega^2 A\tens_A L$ as the image of $1_A$ under the following composition, where the last map is the product in $\Omega^2 B$,
\begin{eqnarray}  \label{cgajhgjcxter}
A \stackrel{\coev} \longrightarrow L^\circ \tens_A L
\stackrel{\id\tens\mathrm{R}_L } \longrightarrow  L^\circ \tens_A \Omega^2 A\tens_A L
\longrightarrow \Omega^2 A\ .
\end{eqnarray}
If we use the curvature in Corollary~\ref{dxfghxfhj} we find $\omega_L$ to be
\begin{eqnarray*}
1_A\longmapsto c\tens c'\longmapsto c\tens c' .\big( 
\extd\eta-\eta\wedge\eta\big).h\tens h'  \longmapsto
\extd\eta-\eta\wedge\eta \ ,
\end{eqnarray*}
where we set $\eta=\zeta-\kappa$. There is no indication that this is a reasonable Chern class. 
However there are several things that we take for granted in ordinary geometry that have to be separately specified in noncommutative geometry. For example, we expect that a left covariant derivative will also have reasonable behaviour for right multiplication -- this is the bimodule covariant derivative condition (see Definition~\ref{ppll}). We have a particularly simple choice of $\sigma$, called $\sigma_0$, in (\ref{zetmghjchdt}).
Corollary~\ref{dxfghjjmxfhj} shows that to get $\sigma=\sigma_0$ in general we need $a.\eta=\eta.a$ for all $a\in A$. Another thing that we take for granted is that the curvature $\mathrm{R}_L$ is a right module map. From Corollary~\ref{dxfghxfhj} this is true when, for all $a\in A$,
$a.(\extd\eta-\eta\wedge\eta)=(\extd\eta-\eta\wedge\eta).a$. If we assume both of these conditions, we deduce that $a.\extd\eta=\extd\eta.a$, and then applying $\extd$ to $a.\eta=\eta.a$ gives $\extd a\wedge\eta+\eta\wedge\extd a=0$ for all $a\in A$.

Subsequently, we deduce that $\eta$ anticommutes with all 1-forms, and in particular that $\eta\wedge\eta=0$. The result is that $\omega_L=\extd\eta$, so $\omega_L$ is in the kernel of $\extd$. However we cannot assume that $[\extd\eta]=0\in H^2_{\mathrm{deRham}}(A)$, as in general $\eta\notin\Omega^1 A$. 
On the other hand we can see that
$[\omega_L]\in H^2_{\mathrm{deRham}}(A)$ is independent of the choice of covariant derivative. Remembering that $\eta=\zeta-\kappa$, where $\zeta\in \Omega^1 A$, the difference of two $\omega_L$ is $\extd$ applied to an element of $\Omega^1 A$. 
Of course the assumptions that the curvature is a right module map, that the connection is a bimodule connection, and (most questionably) that $\sigma=\sigma_0$ are open to question. 
However the framework of line modules does allow calculations to be done in a reasonably sensible manner, and we would hope that it would be used on more examples to give a better idea of what is going on. We conclude with a follow up to Example~\ref{vfauxkyytuyt}, in which we can make the assumptions above and get a sensible answer.

\begin{example} \label{ctfuxkklyfsr}
Following Example~\ref{vfauxkyytuyt}, on $\C_q[SL_2]$ we take the 3D calculus of
\cite{worondiff}. In our conventions this has a basis
\[  e^-=d.\extd b-q b.\extd
d,\quad  e^+=q^{-1} a.\extd c-q^{-2}c.\extd a,\quad  e^0=d.\extd a-q
b.\extd c\] 
of left-invariant 1-forms, is spanned by these as a left
module (according to the above) while the right module relations
and exterior derivative are given in these terms by: 
\begin{eqnarray*}
&& e^\pm
\left(\begin{array}{cc}a & b \\c & d\end{array}\right)
=\left(\begin{array}{cc}qa & q^{-1}b \\qc & q^{-1}d\end{array}\right)
e^\pm,\quad  e^0
\left(\begin{array}{cc}a & b \\c & d\end{array}\right)
=\left(\begin{array}{cc}q^2 a & q^{-2}b \\q^2c & q^{-2}d\end{array}\right)
 e^0 \ ,\cr \cr
&& \extd a=a e^0+q b e^+,\quad \extd
b=a e^--q^{-2}b e^0,\quad \extd c=c  e^0+q d e^+,\quad \extd d=c
e^--q^{-2}d e^0
\end{eqnarray*}
For $\C_q[SL_2]$ the natural extension compatible with
the super-Leibniz rule on higher forms and $\extd^2=0$ is:
\[ \extd  e^0=q^3 e^+\wedge e^-,\quad
\extd e^\pm=\mp q^{\pm 2}[2;q^{- 2}]e^\pm\wedge e^0,\quad
e^\pm\wedge e^\pm= e^0\wedge e^0=0\]
\[
q^2 e^+\wedge e^-+ e^-\wedge e^+=0,\quad  e^0\wedge e ^\pm+q^{\pm
4}e^\pm\wedge e^0=0\] where $[n;q]=(1-q^n)/(1-q)$ denotes a
$q$-integer. This means that there are the same dimensions as
classically, including a unique top form $ e^-\wedge e^+\wedge
e^0$. 

We can set
\begin{eqnarray*}
c\tens c'\ =\ d\tens a- q\,b\tens c\ ,\quad c''\tens c'''\ =\ a\tens d- q^{-1}\, c\tens b\ .
\end{eqnarray*}
From this we calculate $\kappa$ in (\ref{vckuacsxs}), using \cite{BegMa4},
\begin{eqnarray*}
\kappa &=& c.\extd c' \ =\  d.\extd a-q\,b.\extd c \cr
&=& d\,a\,e^0+q\,d\,b\,e^+ - q\,b\,c\,e^0 - q^2\,b\,d\,e^+\cr
&=& (d\,a - q\,b\,c)\,e^0 + (q\,d\,b- q^2\,b\,d)\,e^+\ =\ e^0\ .
\end{eqnarray*}
Now $e^0$ commutes with all elements of $A$, so we can set $\zeta=0$ to get
a metric preserving bimodule connection with (from Corollary~\ref{dxfghjjmxfhj})
\begin{eqnarray*}
\sigma(e\tens\eta) &=& e.\eta.c\tens c' \ ,
\end{eqnarray*}
and curvature given by Corollary~\ref{dxfghxfhj},
and in the case of our example where $\kappa=e^0$ we get
the trace of the curvature being $\omega_l=\extd\kappa=q^3 e^+\wedge e^-$, so the de Rham Chern class would be $q^3 \,[e^+\wedge e^-]\in H^2_{\mathrm{deRham}}(A)$.
\end{example}

\end{document}